\documentclass[12pt]{amsart}
\usepackage{fancybox}
\usepackage{shuffle}
\usepackage{amsmath}
\usepackage{amsfonts}
\usepackage{amssymb}
\usepackage[all]{xy}           
\usepackage{bm}
\usepackage{bbding}
\usepackage{txfonts}
\usepackage{amscd}

\usepackage{xspace}
\usepackage[shortlabels]{enumitem}
\usepackage{ifpdf}

\ifpdf
  \usepackage[colorlinks,final,backref=page,hyperindex]{hyperref}
\else
  \usepackage[colorlinks,final,backref=page,hyperindex,hypertex]{hyperref}
\fi
\usepackage{tikz}
\usepackage[active]{srcltx}

\makeatletter

\newtheorem{thm}{Theorem}[section]
\newtheorem{lem}[thm]{Lemma}
\newtheorem{cor}[thm]{Corollary}
\newtheorem{pro}[thm]{Proposition}
\newtheorem{ex}[thm]{Example}
\theoremstyle{definition}
\newtheorem{rmk}[thm]{Remark}
\newtheorem{defi}[thm]{Definition}

\newcommand{\nc}{\newcommand}
\newcommand{\delete}[1]{}

\nc{\mlabel}[1]{\label{#1}}  
\nc{\mcite}[1]{\cite{#1}}  
\nc{\mref}[1]{\ref{#1}}  
\nc{\mbibitem}[1]{\bibitem{#1}} 

\delete{
\nc{\mlabel}[1]{\label{#1}{\hfill \hspace{1cm}{\bf{{\ }\hfill(#1)}}}}
\nc{\mcite}[1]{\cite{#1}{{\em{{\ }(#1)}}}}  
\nc{\mref}[1]{\ref{#1}{{\em{{\ }(#1)}}}}  
\nc{\mbibitem}[1]{\bibitem[\em #1]{#1}} 
}

\newcommand {\emptycomment}[1]{}

\nc{\oprn}{\theta}
\nc{\Oprn}{\Theta}

\nc{\calo}{\mathcal{O}}
\nc{\oop}{$\mathcal{O}$-operator\xspace}
\nc{\oops}{$\mathcal{O}$-operators\xspace}
\nc{\mrho}{{\bm{\varrho}}}
\nc{\emk}{\mathbf{K}}
\nc{\invlim}{\displaystyle{\lim_{\longleftarrow}}\,}
\nc{\ot}{\otimes}

\newcommand{\lon }{\,\rightarrow\,}
\newcommand{\be }{\begin{equation}}
\newcommand{\ee }{\end{equation}}

\newcommand{\g}{\mathfrak g}
\newcommand{\h}{\mathfrak h}
\newcommand{\m}{\mathfrak m}

\newcommand{\huaB}{\mathcal{B}}
\newcommand{\huaS}{\mathcal{S}}

\newcommand{\huaR}{\mathcal{R}}

\newcommand{\huaG}{\mathcal{G}}

\newcommand{\huaV}{\mathcal{V}}

\newcommand{\huaC}{{\mathcal{C}}}

\newcommand{\huaH}{\mathcal{H}}

\newcommand{\huaO}{{\mathcal{O}}}

\newcommand{\huaZ}{\mathcal{Z}}

\newcommand{\frkg}{\mathfrak g}

\newcommand{\frkC}{\mathfrak C}

\newcommand{\frkI}{\mathfrak I}

\newcommand{\frkL}{\mathfrak L}

\newcommand{\frkR}{\mathfrak R}

\newcommand{\frkT}{\mathfrak T}

\newcommand{\frkX}{\mathfrak X}

\newcommand{\half}{\frac{1}{2}}

\newcommand{\Courant}[1]{\left\llbracket  #1\right\rrbracket }

\newcommand{\Id}{{\rm{Id}}}

\newcommand{\br}[1]{   [ \cdot,    \cdot  ]   }

\newcommand{\Hom}{\mathrm{Hom}}

\newcommand{\Ob}{\mathsf{Ob}}

\newcommand{\gl}{\mathfrak {gl}}

\newcommand{\Ker}{\mathrm{Ker}}

\newcommand{\ad}{\mathrm{ad}}
\newcommand{\pr}{\mathrm{pr}}

\newcommand{\de}{\mathrm{d}}

\newcommand{\LYA}{Lie-Yamaguti algebra}
\nc{\CV}{\mathbf{C}}

\begin{document}

\title[Maurer-Cartan characterization, $L_\infty$-algebras, cohomology of relative RB-operators]{Maurer-Cartan characterization, $L_\infty$-algebras, and cohomology of relative Rota-Baxter operators on Lie-Yamaguti algebras}
\author[Jia Zhao]{Jia Zhao}
\address{Jia Zhao, School of Sciences, Nantong University, Nantong 226019, Jiangsu, China}
\email{zhaojia@ntu.edu.cn}


\author[Yu Qiao]{Yu Qiao*}
\address{Yu Qiao (corresponding author), School of Mathematics and Statistics, Shaanxi Normal University, Xi'an 710119, Shaanxi, China}
\email{yqiao@snnu.edu.cn}

\date{\today}

\begin{abstract}
In this paper, we first construct a {\em differential graded Lie algebra} that controls deformations of a Lie-Yamaguti algebra.
Furthermore, a relative Rota-Baxter operator on a Lie-Yamaguti algebra is characterized as a Maurer-Cartan element
in an appropriate {\em $L_\infty$-algebra} that we build through the graded Lie bracket of \LYA's controlling algebra,
and gives rise to a twisted $L_\infty$-algebra that controls its deformation.
Next we establish the {\em cohomology} theory of relative Rota-Baxter operators on Lie-Yamaguti algebras via the Yamaguti cohomology.
Then we clarify the {\em relationship} between the twisted $L_\infty$-algebra and the cohomology theory.
Finally as byproducts, we classify certain deformations on Lie-Yamaguti algebras using the cohomology theory.
\end{abstract}

\thanks{{\em Mathematics Subject Classification} (2020): 17A30,17B37,17B56}

\keywords{Lie-Yamaguti algebra, Maurer-Cartan element, relative Rota-Baxter operator, $L_\infty$-algebra, cohomology, deformation}

\maketitle

\vspace{-1.1cm}

\tableofcontents

\allowdisplaybreaks

 \section{Introduction}
In this paper, we determine a differential graded Lie algebra as controlling algebra for \LYA s and an $L_\infty$-algebra for their relative Rota-Baxter operators. Then we establish the cohomology theory to classify deformations of relative Rota-Baxter operators.

\subsection{Maurer-Cartan characterizations and deformations}
In mathematics, informally speaking, a deformation of an object is another object that shares the same structure of the original object after a perturbation. Motivated by the foundational work of Kodaira and Spencer \cite{Kodaira} for complex analytic structures, the generalization in algebraic geometry of deformation theory was founded \cite{Hart}. Recently, a number of works referring to deformation quantization turn up in the context of mathematical physics \cite{Kont1,Kont2}. As an application in algebra, Gerstenhaber first studied the deformation theory on associative algebras \cite{Gerstenhaber1,Gerstenhaber2,Gerstenhaber3,Gerstenhaber4,Gerstenhaber5}. Then Nijenhuis and Richardson extended this idea and established the similar results on Lie algebras \cite{Nij1,Nij2}. Deformations of other algebraic structures such as pre-Lie algebras have also been developed \cite{Burde}. In general, deformation theory was set up for binary quadratic operads by Balavoine \cite{bala}.

As we all know, for an algebraic object, deformation theory and cohomology theory are connected tightly. More precisely, a suitable deformation of an object should obey the following rule: on one hand, there exists a differential graded Lie algebra (or generally an $L_\infty$-algebra), to whose Maurer-Cartan elements the algebraic structures correspond; on the other hand, there should be a suitable cohomology theory that characterizes its deformations. Deligne, Drinfeld, and Kontsevich proposed a celebrated slogan that {\em every reasonable deformation is controlled by a differential graded Lie algebra, up to quasi-isomorphisms}. More precisely, this slogan means that for an algebraic structure, Maurer-Cartan elements in a suitable graded Lie algebra are used to characterize realizations of its algebraic structure on a vector space. And moreover, for a given realization of the algebraic structure, a Maurer-Cartan element induces a differential making the graded Lie algebra into a differential graded Lie algebra (this graded Lie algebra is also called the controlling algebra). See \cite{GLST,TBGS,THS,T.S2} for Maurer-Cartan characterizations of some algebraic structures and their relative Rota-Baxter operators.

\subsection{Rota-Baxter operators}
Baxter introduced the notion of Rota-Baxter operators on associative algebras when studying fluctuation theory \cite{Ba}. Then Kupershmidt introduced the notion of $\huaO$-operators (called relative Rota-Baxter operators in the present paper) on Lie algebras when he found that a relative Rota-Baxter operator is a solution to the classical Yang-Baxter equation in \cite{Kupershmidt}, whereas the classical Yang-Baxter equation plays an important role in many fields such as in integration systems \cite{CP}. For more details about the classical Yang-Baxter equation and Rota-Baxter operators, one can see \cite{STS}. Recently, people found that Rota-Baxter operators have many applications. For eaxmple, the first work about Rota-Baxter algebras by Guo had been published \cite{Gub}. From then on, Rota-Baxter algebra open up a new branch in Lie theory and in mathematical physics. Bai and his collaborators examined Rota-Baxter operators in the context of binary operads in \cite{Bai-Bellier-Guo-Ni,PBG}. Moreover, Rota-Baxter operators are closely connected with Hopf algebras.


\subsection{Lie-Yamaguti algebras}
A Lie-Yamaguti algebra is a generalization of a Lie algebra and a Lie triple system, which can be traced back to Nomizu's work on the invariant affine connections on homogeneous spaces in 1950's \cite{Nomizu}. Inspired by thoughts of Nomizu, in 1960's Yamaguti introduced an algebraic object called a general Lie triple system, defined its representation and established its cohomology theory in \cite{Yamaguti1,Yamaguti2}. Kinyon and Weinstein first called this object a \LYA~  when studying Courant algebroids in the earlier 21st century \cite{Weinstein}. Since then, this system is called a \LYA, which has attracted much attention and is widely investigated recently. For instance, Benito, Draper, and Elduque investigated \LYA s related to simple Lie algebras of type $G_2$ \cite{B.D.E}. Afterwards, Benito, Elduque, and Mart$\acute{i}$n-Herce explored irreducible \LYA s in \cite{B.E.M1,B.E.M2}. More recently, Benito, Bremmer, and Madariaga examined orthogonal \LYA s in \cite{B.B.M}. Deformations and extensions of Lie-Yamaguti algebras were explored in \cite{L.CHEN,Zhang1}. Takahashi studied modules over quandles using representations of Lie-Yamaguti algebras in \cite{Takahashi}.

Since a \LYA ~owns two algebraic structures, by the virtue of this object, the first author has been working on the deformation theory of \LYA s these years. For example, Sheng, the first author, and Zhou analyzed product structures and complex structures on Lie-Yamaguti algebras by means of Nijenhuis operators in \cite{Sheng Zhao}.  Afterwards, Sheng and the first author introduced the notion of relative Rota-Baxter operators on Lie-Yamaguti algebras and revealed the fact that a pre-Lie-Yamaguti algebra is the underlying algebraic structure of relative Rota-Baxter operators \cite{SZ1}. Since then, thanks to importance of deformations and relative Rota-Bxter operators on \LYA s, we established Lie-Yamaguti bialgebra theory and clarified the relationship between the solution to the classical Lie-Yamaguti Yang-Baxter equation and relative Rota-Baxter operators in \cite{ZQ1}, and studied cohomology and deformations of $\mathsf{LieYRep}$ pairs and then explored several properties of relative Rota-Baxter-Nijenhuis structures in \cite{ZQ2}.

{\bf Motivation.} Since Sheng and the first author defined relative Rota-Baxter operators on \LYA s \cite{SZ1}, according to Deligne, Drinfeld, and Kontsevich's slogan, it is natural to consider Maurer-Cartan characterizations for relative Rota-Baxter operators on \LYA s and its relationship between cohomology and the Maurer-Cartan elements. However, the first and foremost, we have to consider Maurer-Cartan characterization for \LYA s because on the one hand, unlike associative algebras, Lie algebras, Leibniz algebras, or even $3$-Lie algebras, no works of controlling algebras of \LYA s appeared before; on the other hand, it is inevitable, for we construct the controlling algebra for relative Rota-Baxter operators. Thus, this step of constructing controlling algebra for \LYA s also promotes development of deformation theory for \LYA s. Anyway, our project can be transferred to the following questions:
\begin{itemize}
\item Does there exist a suitable graded Lie algebra whose Maurer-Cartan elements correspond to the \LYA ~structures on a vector space?
\item Does there exist a suitable algebra whose Maurer-Catan elemnts are precisely the relative Rota-Baxter operators on \LYA s?
\item Does there exist an appropriate cohomology theory of relative Rota-Baxter operators on Lie-Yamaguti algebras, which can be used to classify certain types of deformations?
\end{itemize}
After solving these problems, we are able to consider deformations spontaneously, perfecting researches about cohomology and deformation theory on \LYA s and on their relative Rota-Baxter operators.

\subsection{Outline of the paper}
We tackle these problems as follows. After some preliminaries introduced in Section 2, Section 3 answers the first question. As was stated before, no works of controlling algebras of \LYA s appeared, thus constructing a graded Lie algebra whose Maurer-Cartan elements correspond to the \LYA ~structures is a {\em difficulty} to overcome. Once the graded Lie algebra is constructed (see Theorem \ref{fund}), consequently a Maurer-Cartan element induces a differential making the graded Lie algebra into a differential graded Lie algebra as the controlling algebra for \LYA s (see Theorem \ref{fund2}). Moreover, we establish the {\em relation} between the controlling algebra and the coboundary operator associated to the adjoint representation (see Theorem \ref{diffe}).

Section 4 is the solution to the last two questions, which is the core of the present paper. Immediately after the constructing controlling algebra for \LYA s, we make use of the graded Lie brackets by the method of derived brackets to construct an $L_\infty$-algebra whose Maurer-Cartan elements are precisely relative Rota-Baxter operators (see Theorem \ref{main}). Moreover, from a relative Rota-Baxter operator (as a Maurer-Cartan element), we obtain a twisted $L_\infty$-algebra that controls deformations of relative Rota-Baxter operators (see Theorem \ref{deformation}). As an application, we endow a continuous function space and a series ring with \LYA ~structure respectively and give two integral operators which are treated as relative Rota-Baxter operators on them. And then we build cohomology of relative Rota-Baxter operators as follows. Given a relative Rota-Baxter operator $T:V\longrightarrow\g$ with respect to a representation $(V;\rho,\mu)$, we mentioned in \cite{SZ1} that there is a \LYA ~structure $([\cdot,\cdot]_T,\Courant{\cdot,\cdot,\cdot}_T)$ on $V$, in which case the \LYA ~$(V,[\cdot,\cdot]_T,\Courant{\cdot,\cdot,\cdot}_T)$ is called the sub-adjacent \LYA. In order to establish the cohomology theory, we first should build the cohomology of this sub-adjacent \LYA ~using Yamaguti cohomology. Thus, we have to construct a representation of the sub-adjacent \LYA ~(see Theorem \ref{represent}). 
However, note that the cochain complex of Yamaguti cohomology starts only from $1$-cochians, {\em not} from $0$-cochains. We have to point out that $0$-cochains play an import role in our study of deformations.
Thus the next step is to choose $0$-cochains appropriately and build a proper coboundary map from the set of $0$-cochains to that of $1$-cochains, which is another {\em difficulty} (see Proposition \ref{0cocy}).
In this way, we obtain a cochain complex (associated to $V$) starting from $0$-cochains,
which gives rise to the cohomology of the relative Rota-Baxter operator $T$ on Lie-Yamaguti algebras $(\g,[\cdot,\cdot],\Courant{\cdot,\cdot,\cdot})$ (see Definition \ref{cohomology}). We also establish the {\em relation} between the differential of the twisted $L_\infty$-algebra and the coboundary operator of cohomology of relative Rota-Baxter operators parallel to the case of \LYA s (see Theorem \ref{diff}).

Finally in Section 5, as an application, we make use of the cohomology theory to investigate deformations of relative Rota-Baxter operators on Lie-Yamaguti algebras. We intend to consider two kinds of deformations: linear and higher order deformations. It turns out that our cohomology theory satisfies the rule that is mentioned above and works well (see Theorem \ref{thm1}  and Theorem \ref{ob}).


A Lie-Yamaguti algebra owns two algebraic operations, which makes the operation of its controlling algebra and cochain complex much more complicated than other algebras, such as Lie algebras, pre-Lie algebras, Leibniz algebras or even $3$-Lie algebras, one of which owns only one operation. As a result, the computation is technically difficult in constructing the controlling algebras and defining the cohomology of relative Rota-Baxter operators. Note that a Lie triple system is a spacial case of a Lie-Yamaguti algebra \cite{Lister}, so the conclusions in the present paper can also be  adapted to the Lie triple system context. 

\emptycomment{
The paper is structured as follows.
In Section 2, we recall some basic concepts including those of Lie-Yamaguti algebras, representations and cohomology.
In Section 3, the cohomology theory of relative Rota-Baxter operators on Lie-Yamaguti algebras is constructed by using that of Lie-Yamaguti algebras.
Finally, we utilize our established cohomology theory to analyze three kinds of deformations of relative Rota-Baxter operators on Lie-Yamaguti algebras: namely linear, formal, and higher order deformations in Section 4.}

 In this paper, all the vector spaces are over $\mathbb{K}$, a field of characteristic $0$.

\section{Preliminaries}
In this section, we recall some basic notions such as Lie-Yamaguti algebras, representations and their cohomology theory.
The notion of  Lie-Yamaguti algebras was first defined by  Yamaguti in \cite{Yamaguti1,Yamaguti2}.
\begin{defi}\cite{Weinstein}\label{LY}
A {\bf Lie-Yamaguti algebra} is a vector space $\g$, together with a bilinear bracket $[\cdot,\cdot]:\wedge^2  \mathfrak{g} \to \mathfrak{g} $ and a trilinear bracket $\Courant{\cdot,\cdot,\cdot}:\wedge^2\g \otimes  \mathfrak{g} \to \mathfrak{g} $ such that the following equations are satisfied for all $x,y,z,w,t \in \g$,
\begin{eqnarray}
~ &&\label{LY1}[[x,y],z]+[[y,z],x]+[[z,x],y]+\Courant{x,y,z}+\Courant{y,z,x}+\Courant{z,x,y}=0,\\
~ &&\Courant{[x,y],z,w}+\Courant{[y,z],x,w}+\Courant{[z,x],y,w}=0,\\
~ &&\label{LY3}\Courant{x,y,[z,w]}=[\Courant{x,y,z},w]+[z,\Courant{x,y,w}],\\
~ &&\Courant{x,y,\Courant{z,w,t}}=\Courant{\Courant{x,y,z},w,t}+\Courant{z,\Courant{x,y,w},t}+\Courant{z,w,\Courant{x,y,t}}.\label{fundamental}
\end{eqnarray}
\end{defi}

\begin{ex}\cite{Nomizu}
Let $M$ be a closed manifold 
with an affine connection, and denote by $\frkX(M)$ the set of vector fields on $M$. For all $x,y,z\in \frkX(M)$, set
\begin{eqnarray*}
[x,y]=-T(x,y),\quad \Courant{x,y,z}=-R(x,y)z,
\end{eqnarray*}
where $T$ and $R$ are torsion tensor and curvature tensor respectively. It turns out that the triple
$ (\frkX(M),[\cdot,\cdot],\Courant{\cdot,\cdot,\cdot})$ forms an (infinite-dimensional) \LYA.
\end{ex}

\begin{defi}\cite{Yamaguti2}\label{defi:representation}
Let $(\g,[\cdot,\cdot],\Courant{\cdot,\cdot,\cdot})$ be a Lie-Yamaguti algebra and $V$ a vector space. A {\bf representation}  of $\g$ on $V$ consists of a linear map $\rho:\g \to \gl(V)$ and a bilinear map $\mu:\otimes^2 \g \to \gl(V)$ such that for all $x,y,z,w \in \g$,
\begin{eqnarray}
~&&\label{RLYb}\mu([x,y],z)-\mu(x,z)\rho(y)+\mu(y,z)\rho(x)=0,\\
~&&\label{RLYd}\mu(x,[y,z])-\rho(y)\mu(x,z)+\rho(z)\mu(x,y)=0,\\
~&&\label{RLYe}\rho(\Courant{x,y,z})=[D_{\rho,\mu}(x,y),\rho(z)],\\
~&&\label{RYT4}\mu(z,w)\mu(x,y)-\mu(y,w)\mu(x,z)-\mu(x,\Courant{y,z,w})+D_{\rho,\mu}(y,z)\mu(x,w)=0,\\
~&&\label{RLY5}\mu(\Courant{x,y,z},w)+\mu(z,\Courant{x,y,w})=[D_{\rho,\mu}(x,y),\mu(z,w)],
\end{eqnarray}
where the bilinear map $D_{\rho,\mu}:\otimes^2\g \to \gl(V)$ is given by
\begin{eqnarray*}
 D_{\rho,\mu}(x,y):=\mu(y,x)-\mu(x,y)+[\rho(x),\rho(y)]-\rho([x,y]), \quad \forall x,y \in \g.\label{rep}
 \end{eqnarray*}
It is obvious that $D_{\rho,\mu}$ is skew-symmetric. We write $D$ in the sequel without ambiguities and we denote a representation of $\g$ on $V$ by $(V;\rho,\mu)$.
\end{defi}
\emptycomment{
\begin{rmk}\label{rmk:rep}
Let $(\g,[\cdot,\cdot],\Courant{\cdot,\cdot,\cdot})$ be a Lie-Yamaguti algebra and $(V;\rho,\mu)$ its representation. If $\rho=0$ and the Lie-Yamaguti algebra $\g$ reduces to a Lie tripe system $(\g,\Courant{\cdot,\cdot,\cdot})$,  then $(V;\mu)$  is a representation of the Lie triple systems $(\g,\Courant{\cdot,\cdot,\cdot})$; If $\mu=0$, $D=0$ and the Lie-Yamaguti algebra $\g$ reduces to a Lie algebra $(\g,[\cdot,\cdot])$, then $(V;\rho)$ is a representation  of the Lie algebra $(\g,[\cdot,\cdot])$. So the above definition of a representation of a Lie-Yamaguti algebra is a natural generalization of representations of Lie algebras and Lie triple systems.
\end{rmk}}

By a direct computation, we have

\begin{pro}
If $(V;\rho,\mu)$ is a representation of a Lie-Yamaguti algebra $(\g,[\cdot,\cdot],\Courant{\cdot,\cdot,\cdot})$, then we have the following equalities:
\begin{eqnarray}
\label{RLYc}&&D([x,y],z)+D([y,z],x)+D([z,x],y)=0;\\
\label{RLY5a}&&D(\Courant{x,y,z},w)+D(z,\Courant{x,y,w})=[D(x,y),D(z,w)];\\
~ &&\mu(\Courant{x,y,z},w)=\mu(x,w)\mu(z,y)-\mu(y,w)\mu(z,x)-\mu(z,w)D(x,y).\label{RLY6}
\end{eqnarray}
\end{pro}

\begin{ex}\label{ad}
Let $(\g,[\cdot,\cdot],\Courant{\cdot,\cdot,\cdot})$ be a Lie-Yamaguti algebra. We define linear maps $\ad:\g \to \gl(\g)$ and $\frkR :\otimes^2\g \to \gl(\g)$ to be $x \mapsto \ad_x$ and $(x,y) \mapsto \mathfrak{R}_{x,y}$, where $\ad_xz=[x,z]$ and $\mathfrak{R}_{x,y}z=\Courant{z,x,y}$ for all $z \in \g$ respectively. Then $(\g;\ad,\mathfrak{R})$ forms a representation of $\g$ on itself, called the {\bf adjoint representation}. In this case, $\frkL\triangleq D_{\ad,\frkR}$ is given by for all $x,y \in \g$,
\begin{eqnarray*}
\frkL_{x,y}=\mathfrak{R}_{y,x}-\mathfrak{R}_{x,y}+[\ad_x,\ad_y]-\ad_{[x,y]}.
\end{eqnarray*}
By \eqref{LY1}, we have
\begin{eqnarray*}
\frkL_{x,y}z=\Courant{x,y,z}, \quad \forall z \in \g.\label{lef}
\end{eqnarray*}
\end{ex}

Representations of a Lie-Yamaguti algebra can be characterized by the semidirect product Lie-Yamaguti algebras.

\begin{pro}\label{semi}\cite{Zhang1}
Let $(\g,[\cdot,\cdot],\Courant{\cdot,\cdot,\cdot})$ be a Lie-Yamaguti algebra and $V$ a vector space. Let $\rho:\g \to \gl(V)$ and $\mu:\otimes^2 \g \to \gl(V)$ be linear maps. Then $(V;\rho,\mu)$ is a representation of $(\g,[\cdot,\cdot],\Courant{\cdot,\cdot,\cdot})$ if and only if there is a Lie-Yamaguti algebra structure $([\cdot,\cdot]_{\rho,\mu},\Courant{\cdot,\cdot,\cdot}_{\rho,\mu})$ on the direct sum $\g \oplus V$ which is defined by for all $x,y,z \in \g, ~u,v,w \in V$,
\begin{eqnarray*}
\label{semi1}[x+u,y+v]_{\rho,\mu}&=&[x,y]+\rho(x)v-\rho(y)u,\\
\label{semi2}~\Courant{x+u,y+v,z+w}_{\rho,\mu}&=&\Courant{x,y,z}+D(x,y)w+\mu(y,z)u-\mu(x,z)v,
\end{eqnarray*}
This Lie-Yamaguti algebra $(\g \oplus V,[\cdot,\cdot]_{\rho,\mu},\Courant{\cdot,\cdot,\cdot}_{\rho,\mu})$ is called the {\bf semidirect product Lie-Yamaguti algebra}, and denoted by $\g \ltimes_{\rho,\mu} V$.
\end{pro}

The cohomology theory of Lie-Yamaguti algebras was established in \cite{Yamaguti2}. Let $(\g,[\cdot,\cdot],\Courant{\cdot,\cdot,\cdot})$ be a \LYA ~and $(V;\rho,\mu)$ a representation. We denote by $C^p_{\rm LieY}(\g,V)~(p \geqslant 1)$ the set of $p$-cochains, where
\begin{eqnarray*}
C^{n+1}_{\rm LieY}(\g,V)\triangleq
\begin{cases}
\Hom(\underbrace{\wedge^2\g\otimes \cdots \otimes \wedge^2\g}_n,V)\times \Hom(\underbrace{\wedge^2\g\otimes\cdots\otimes\wedge^2\g}_{n}\otimes\g,V), & n\geqslant 1,\\
\Hom(\g,V), &n=0.
\end{cases}
\end{eqnarray*}

For $p\geqslant 1$, the coboundary operator $\delta:C^p_{\rm LieY}(\g,V)\to C^{p+1}_{\rm LieY}(\g,V)$ is defined as follows:
\begin{itemize}
\item If $n\geqslant 1$, for any $F=(f,g)\in C^{n+1}_{\rm LieY}(\g,V)$, the coboundary map
$$\delta=(\delta_{\rm I},\delta_{\rm II}):C^{n+1}_{\rm LieY}(\g,V)\to C^{n+2}_{\rm LieY}(\g,V),$$
$$\qquad \qquad\qquad \qquad\qquad \quad F\mapsto(\delta_{\rm I}(F),\delta_{\rm II}(F)),$$
 is given by:
\begin{eqnarray}
~ &&\nonumber\Big(\delta_{\rm I}(F)\Big)(\frkX_1,\cdots,\frkX_{n+1})\\
~\label{cohomolo1} &=&(-1)^n\Big(\rho(x_{n+1})g(\frkX_1,\cdots,\frkX_n,y_{n+1})-\rho(y_{n+1})g(\frkX_1,\cdots,\frkX_n,x_{n+1})\\
~ &&\nonumber-g(\frkX_1,\cdots,\frkX_n,[x_{n+1},y_{n+1}])\Big)\\
~ &&\nonumber+\sum_{k=1}^{n}(-1)^{k+1}D(\frkX_k)f(\frkX_1,\cdots,\hat{\frkX_k},\cdots,\frkX_{n+1})\\
~ &&\nonumber+\sum_{1\leqslant k<l\leqslant n+1}(-1)^{k}f(\frkX_1,\cdots,\hat{\frkX_k},\cdots,\frkX_k\circ\frkX_l,\cdots,\frkX_{n+1}),
\end{eqnarray}
and
\begin{eqnarray}
~ &&\nonumber\Big(\delta_{\rm II}(F)\Big)(\frkX_1,\cdots,\frkX_{n+1},z)\\
~ \label{cohomolo2}&=&(-1)^n\Big(\mu(y_{n+1},z)g(\frkX_1,\cdots,\frkX_n,x_{n+1})-\mu(x_{n+1},z)g(\frkX_1,\cdots,\frkX_n,y_{n+1})\Big)\\
~ &&\nonumber+\sum_{k=1}^{n+1}(-1)^{k+1}D(\frkX_k)g(\frkX_1,\cdots,\hat{\frkX_k},\cdots,\frkX_{n+1},z)\\
~ &&\nonumber+\sum_{1\leqslant k<l\leqslant n+1}(-1)^kg(\frkX_1,\cdots,\hat{\frkX_k},\cdots,\frkX_k\circ\frkX_l,\cdots,\frkX_{n+1},z)\\
~ &&\nonumber+\sum_{k=1}^{n+1}(-1)^kg(\frkX_1,\cdots,\hat{\frkX_k},\cdots,\frkX_{n+1},\Courant{x_k,y_k,z}),
\end{eqnarray}
where $\frkX_i=x_i\wedge y_i\in\wedge^2\g,~(i=1,\cdots,n+1),~z\in \g$ and the operation $\circ$ means that
 $$\frkX_k\circ\frkX_l\triangleq\Courant{x_k,y_k,x_l}\wedge y_l+x_l\wedge\Courant{x_k,y_k,y_l}.$$

\item For the case that $n=0$, any element $f \in C^1_{\rm LieY}(\g,V)$ given, the coboundary map
$$\delta:C^1_{\rm LieY}(\g,V)\to C^2_{\rm LieY}(\g,V),$$
$$\qquad \qquad \qquad f\mapsto (\delta_{\rm I}(f),\delta_{\rm II}(f)),$$
is given by
\begin{eqnarray}
\label{1cochain}(\delta_{\rm I}(f))(x,y)&=&\rho(x)f(y)-\rho(y)f(x)-f([x,y]),\\
~ \label{2cochain}(\delta_{\rm II}(f))(x,y,z)&=&D(x,y)f(z)+\mu(y,z)f(x)-\mu(x,z)f(y)-f(\Courant{x,y,z}),\quad \forall x,y, z\in \g.
\end{eqnarray}
\end{itemize}

Yamaguti showed in \cite{Yamaguti2} that $\delta$ is a differential, i.e., $\delta\circ\delta=0$. More precisely, for any $f\in C^1_{\rm LieY}(\g,V)$, we have
$$\delta_{\rm I}\Big(\delta_{\rm I}(f),\delta_{\rm II}(f)\Big)=0\quad\text{and}\quad\delta_{\rm II}\Big(\delta_{\rm I}(f),\delta_{\rm II}(f)\Big)=0.$$
Moreover, for all $F\in C^p_{\rm LieY}(\g,V)~(p\geqslant 2)$, we have
$$\delta_{\rm I}\Big(\delta_{\rm I}(F),\delta_{\rm II}(F)\Big)=0\quad\text{and}\quad\delta_{\rm II}\Big(\delta_{\rm I}(F),\delta_{\rm II}(F)\Big) =0.$$
Thus the cochain complex $(C^\bullet_{\rm LieY}(\g,V)=\bigoplus_{p=1}^\infty C^p_{\rm LieY}(\g,V),\delta)$ is well defined, whose cohomology is called the {\bf Yamaguti cohomology} in this paper.

\begin{defi}
An $p$-cocycle $(p\geqslant 1)$ is an element $F=(f,g)$ in $C^p_{\rm LieY}(\g,V)~(p\geqslant 2)$ (resp. $f\in C^1_{\rm LieY}(\g,V)$) such that $\delta(F)=0$ (resp. $\delta(f)=0$). The set of $p$-cocycles is denoted by $Z^p_{\rm LieY}(\g,V)$; for an element $F$ in $C^p_{\rm LieY}(\g,V)~(p\geqslant 2)$, if there exists $G=(h,s)\in C^{p-1}_{\rm LieY}(\g,V)$~(resp. $t\in C^1(\g,V)$, if $p=2$) such that $F=\delta(G)$~(resp. $F=\delta(t)$), then $F$ is called an $p$-coboundary. The set of $p$-coboundaries is denoted by $B^p_{\rm LieY}(\g,V)$. The resulting $p$-cohomology group is defined by the factor space
$$H^p_{\rm LieY}(\g,V)=Z^p_{\rm LieY}(\g,V)/B^p_{\rm LieY}(\g,V)~\quad (p\geqslant 2).$$ 
\end{defi}

\section{Maurer-Cartan characterizations for Lie-Yamaguti algebras}


In this section, for a preparation of constructing controlling algebra for relative Rota-Baxter operators, we have to construct the controlling algebra for Lie-Yamaguti algebras first. Let us recall some notions and basic facts in \cite{Loday}.

A degree $1$ element $x\in \g_1$ is called a {\bf Maurer-Cartan} element of a differential graded Lie algebra $(\g=\oplus_{k\in\mathbb Z}\g_k,[\cdot,\cdot],d)$ if it satisfies the Maurer-Cartan equation: $dx+\half [x,x]=0.$ Note that a graded Lie algebra is a special differential graded Lie algebra with $d=0$. Correspondingly, for a graded Lie algebra $(\g=\oplus_{k\in\mathbb Z}\g_k,[\cdot,\cdot])$, an element $x\in \g_1$ satisfying $[x,x]=0$ is a Maurer-Cartan element of $\g$.

A permutation $\sigma\in S_n$ is called an {\bf $(i,n-i)$-shuffle} if $\sigma(1)<\cdots <\sigma(i)$ and $\sigma(i+1)<\cdots <\sigma(n)$. If $i=0$ or $i=n$, we assume $\sigma=\Id$. The set of $(i,n-i)$-shuffles is denoted by $\mathbb S_{(i,n-i)}$.

Let $\g$ be a vector space. Denote by
\begin{eqnarray*}
\frkC^p(\g,\g):=
\begin{cases}
\Hom(\underbrace{(\wedge^2\g)\otimes\cdots(\wedge^2\g)}_p,\g)\oplus\Hom(\underbrace{(\wedge^2\g)\otimes\cdots(\wedge^2\g)}_p\otimes\g,\g),& p\geqslant 1,\\
\Hom(\g,\g),&p=0.
\end{cases}
\end{eqnarray*}
Then $\frkC^\bullet(\g,\g)=\oplus_{p\geqslant 0}\frkC^p(\g,\g)$ is a graded vector space, in which the degree of elements in $\frkC^p(\g,\g)$ is assumed to be $p$. For $P=(P_{\rm I},P_{\rm II}) \in \frkC^p(\g,\g),~Q=(Q_{\rm I},Q_{\rm II})\in \frkC^q(\g,\g)~(p,q\geqslant 1)$, we define $P\circ Q=\Big((P\circ Q)_{\rm I},(P\circ Q)_{\rm II}\Big)\in \frkC^{p+q}(\g,\g)$ to be
{\footnotesize
\begin{eqnarray}
~ \nonumber&&\Big(P\circ Q\Big)_{\rm I}(\frkX_1,\cdots,\frkX_{p+q})\\
~ \nonumber&=&\sum_{\sigma\in\mathbb S_{(p,q)}\atop \sigma(p+q)=p+q}(-1)^{pq}sign(\sigma)P_{\rm II}(\frkX_{\sigma(1)},\cdots,\frkX_{\sigma(p)},Q_{\rm I}(\frkX_{\sigma(p+1)},
\cdots,\frkX_{\sigma(p+q)}))\\
 \label{gradedbra1}&&+\sum_{k=1}^p(-1)^{(k-1)q}\sum_{\sigma\in \mathbb S_{(k-1,q)}}sign(\sigma)P_{\rm I}(\frkX_{\sigma(1)},\cdots,\frkX_{\sigma(k-1)},x_{q+k}\wedge
Q_{\rm II}(\frkX_{\sigma(k)},\cdots,\frkX_{\sigma(k+q-1)},y_{k+q}),\frkX_{k+q+1},\cdots,\frkX_{p+q})\\
~ \nonumber&&+\sum_{k=1}^p(-1)^{(k-1)q}\sum_{\sigma\in \mathbb S_{(k-1,q)}}sign(\sigma)P_{\rm I}(\frkX_{\sigma(1)},\cdots,\frkX_{\sigma(k-1)},
Q_{\rm II}(\frkX_{\sigma(k)},\cdots,\frkX_{\sigma(k+q-1)},x_{k+q})\wedge y_{k+q},\frkX_{k+q+1},\cdots,\frkX_{p+q}),
\end{eqnarray}}
and
{\footnotesize
\begin{eqnarray}
~ \nonumber&&\Big(P\circ Q\Big)_{\rm II}(\frkX_1,\cdots,\frkX_{p+q},x)\\
~ \nonumber&=&\sum_{\sigma\in\mathbb S_{(p,q)}}(-1)^{pq}sign(\sigma)P_{\rm II}(\frkX_{\sigma(1)},\cdots,\frkX_{\sigma(p)},Q_{\rm II}(\frkX_{\sigma(p+1)},
\cdots,\frkX_{\sigma(p+q)},x))\\
\label{gradedbra3}&&+\sum_{k=1}^p(-1)^{(k-1)q}\sum_{\sigma\in \mathbb S_{(k-1,q)}}sign(\sigma)P_{\rm II}(\frkX_{\sigma(1)},\cdots,\frkX_{\sigma(k-1)},x_{q+k}\wedge
Q_{\rm II}(\frkX_{\sigma(k)},\cdots,\frkX_{\sigma(k+q-1)},y_{k+q}),\frkX_{k+q+1},\cdots,\frkX_{p+q},x)\\
~ \nonumber&&+\sum_{k=1}^p(-1)^{(k-1)q}\sum_{\sigma\in \mathbb S_{(k-1,q)}}sign(\sigma)P_{\rm II}(\frkX_{\sigma(1)},\cdots,\frkX_{\sigma(k-1)},
Q_{\rm II}(\frkX_{\sigma(k)},\cdots,\frkX_{\sigma(k+q-1)},x_{k+q})\wedge y_{k+q},\frkX_{k+q+1},\cdots,\frkX_{p+q},x),
\end{eqnarray}}
 In particular, for $f,g\in \frkC^0(\g,\g)=\Hom(\g,\g)$ and $P=(P_{\rm I},P_{\rm II})\in \frkC^p(\g,\g)~(p\geqslant 1)$, we define
 \begin{eqnarray}
 \label{gradedbra4}\Big(P\circ f\Big)_{\rm I}(\frkX_1,\cdots,\frkX_p)&=&\sum_{k=1}^pP_{\rm I}\Big(\frkX_1,\cdots,\frkX_{k-1},x_k\wedge f(y_k),\frkX_{k+1},\cdots,\frkX_p\Big)\\
 ~\nonumber&& +\sum_{k=1}^pP_{\rm I}\Big(\frkX_1,\cdots,\frkX_{k-1},f(x_k)\wedge y_k,\frkX_{k+1},\cdots,\frkX_p\Big),\\
 ~ \nonumber&&\\
 \Big(f\circ P\Big)_{\rm I}(\frkX_1,\cdots,\frkX_p)&=&f\Big(P_{\rm I}(\frkX_1,\cdots,\frkX_p)\Big),\label{gradedbra5}\\
 \Big(P\circ f\Big)_{\rm II}(\frkX_1,\cdots,\frkX_p,x)&=&\sum_{k=1}^pP_{\rm II}\Big(\frkX_1,\cdots,\frkX_{k-1},x_k\wedge f(y_k),\frkX_{k+1},\cdots,\frkX_{p},x\Big)\label{gradedbra6}\\
 ~\nonumber&& +\sum_{k=1}^pP_{\rm II}\Big(\frkX_1,\cdots,\frkX_{k-1},f(x_k)\wedge y_k,\frkX_{k+1},\cdots,\frkX_{p},x\Big)\\
 ~ \nonumber&&+P_{\rm II}\Big(\frkX_1,\cdots,\frkX_p,f(x)\Big),\\
 \label{gradedbra2}\Big(f\circ P\Big)_{\rm II}(\frkX_1,\cdots,\frkX_p,x)&=&f\Big(P_{\rm II}(\frkX_1,\cdots,\frkX_p,x)\Big).
 \end{eqnarray}
 Moreover, $f\circ g$ means the composition of $f$ and $g$.

Let us introduce some notations. Let $\g$ and $V$ be vector spaces. Denote a degree $1$ element $(\pi,\omega)\in \frkC^1(\g,\g)$ in $\frkC^\bullet(\g,\g)$ by $\Pi$. Correspondingly, for another element $\Pi'=(\pi',\omega')\in \frkC^1(\g,\g)$, the element $(\pi+\pi',\omega+\omega')$ is written as $\Pi+\Pi'$. Moreover, recall that
$$\frkC^1(\g\oplus V,\g\oplus V)=\Hom(\wedge^2(\g\oplus V),\g \oplus V)\oplus \Hom(\wedge^2(\g\oplus V)\otimes (\g\oplus V),\g\oplus V).$$
For linear maps $\rho:\g\to\gl(V)$ and $\mu:\otimes^2\g\to \gl(V)$, define the degree $1$ elements $\overline\Pi=(\bar\pi,\bar\omega)$ and $\overline\Theta=(\bar\rho,\bar\mu)\in \frkC^1(\g\oplus V,\g \oplus V)$ to be for all $x,y,z\in \g$ and $u,v,w\in V$,
\begin{eqnarray*}
\begin{cases}
\bar{\pi}(x+u,y+v)=\pi(x,y),\\
\bar{\omega}(x+u,y+v,z+w)=\omega(x,y,z),
\end{cases}
\end{eqnarray*}
and
\begin{eqnarray*}
\begin{cases}
\bar{\rho}(x+u,y+v)=\rho(x)v-\rho(y)u,\\
\bar{\mu}(x+u,y+v,z+w)=D(x,y)w+\mu(y,z)u-\mu(x,z)v,
\end{cases}
\end{eqnarray*}
respectively. In the following, we denote $\bar \rho$ and $\bar\mu$ by $\rho$ and $\mu$ respectively,
then the degree $1$ element $\overline\Pi+\overline\Theta=(\bar\pi+\bar\rho,\bar\omega+\bar\mu)\in \frkC^1(\g\oplus V,\g\oplus V)$ is written as $\Pi+\Theta=(\pi+\rho,\omega+\mu)$ correspondingly. In the sequel, we also denote a \LYA ~$(\g,[\cdot,\cdot],\Courant{\cdot,\cdot})$ by $(\g,\pi,\omega)$, where $\pi=[\cdot,\cdot]$ and $\omega=\Courant{\cdot,\cdot,\cdot}$.

\begin{thm}\label{fund}
With the above notations, the graded vector space $\frkC^\bullet(\g,\g)$ equipped with the graded commutator
$$[P,Q]_{\mathsf{LieY}}=P\circ Q-(-1)^{pq}Q\circ P$$
is a graded Lie algebra, where the operation $\circ$ is defined as \eqref{gradedbra1}-\eqref{gradedbra2}. Moreover, if
$\Pi\in \frkC^1(\g,\g)$ defines a Lie-Yamaguti algebra structure on $\g$, then $\Pi$ is a Maurer-Cartan element of the graded Lie algebra $(C^\bullet(\g,\g),[\cdot,\cdot]_{\mathsf{LieY}})$.
\end{thm}
\begin{proof}
We first prove the graded bracket is skew-symmetric.
Indeed, for all $P \in \frkC^p(\g,\g),~Q \in \frkC^q(\g,\g)~(p,q\geqslant 0)$, we write
$$[P,Q]_{\mathsf{LieY}}=([P,Q]_{\rm I},[P,Q]_{\rm II})$$
where
\begin{eqnarray*}
~[P,Q]_{\rm I}&=&\Big(P\circ Q\Big)_{\rm I}-(-1)^{pq}\Big(Q\circ P\Big)_{\rm I}\\
~[P,Q]_{\rm II}&=&\Big(P\circ Q\Big)_{\rm II}-(-1)^{pq}\Big(Q\circ P\Big)_{\rm II}.
\end{eqnarray*}
We have
\begin{eqnarray*}
~[P,Q]_{\mathsf{LieY}}&=&([P,Q]_{\rm I},[P,Q]_{\rm II})\\
~ &=&(-(-1)^{pq}[Q,P]_{\rm I},-(-1)^{pq}[Q,P]_{\rm II}\\
~ &=&-(-1)^{pq}\Big([Q,P]_{\rm I},[Q,P]_{\rm II}\Big)\\
~ &=&-(-1)^{pq}[Q,P]_{\mathsf{LieY}}.
\end{eqnarray*}

Next we prove that the bracket $[\cdot,\cdot]_{\mathsf{LieY}}$ satisfies the graded Jacobi Identity, which means that both $[\cdot,\cdot]_{\rm I}$ and $[\cdot,\cdot]_{\rm {II}}$ satisfy the graded Jacobi Identity. In fact, that the bracket $[\cdot,\cdot]_{\rm {II}}$ satisfies the graded Jacobi Identity has been proved in \cite{Rot} using the weight rule. For any element $P\in \frkC^p(\g,\g)~(p\geqslant 1)$, fix a nonzero element $x\in \g$, and define a map
$$\Phi:\Hom(\wedge^2\g^{\otimes p}\otimes\g,\g)\to\Hom(\wedge^2\g^{\otimes p},\g)$$
to be
$$\Phi(P)(\frkX_1,\cdots,\frkX_p)=P(\frkX_1,\cdots,\frkX_p,x),\quad \forall P\in \Hom(\wedge^2\g^{\otimes p}\otimes\g,\g),$$
where $\frkX_k\in \wedge^2\g~(k=1,2,\cdots,p)$. Then since the graded vector space $\oplus_{p\geqslant 0}\Hom(\wedge^2\g^{\otimes p}\otimes\g,\g)$ endowed with $[\cdot,\cdot]_{\rm II}$ is a graded Lie algebra, the graded vector space $\oplus_{p\geqslant 0}\Hom(\wedge^2\g^{\otimes p},\g)$ is also a graded Lie algebra via the map $\Phi$, whose graded Lie bracket is exactly the bracket $[\cdot,\cdot]_{\rm I}$ for two variables whose degrees $\geqslant 1$. Furthermore, by a direct computation, if either two variables is of degree $0$, then the bracket $[\cdot,\dot]_{\rm I}$ is also a graded Lie bracket. Thus the bracket $[\cdot,\cdot]_{\mathsf{LieY}}$ satisfies the graded Jacobi Identity.

Finally, for all $x,y,z,w,t\in \g$, we have
\begin{eqnarray*}
~ &&[\Pi,\Pi]_{\rm I}(x,y,z,w)\\
~ &=&2\Big(\big(\Pi\circ\Pi\big)_{\rm I}(x,y,z,w)\Big)\\
~ &=&2\Big(\omega(x,y,\pi(z,w))-\pi(\omega(x,y,z),w)-\pi(z,\omega(x,y,w))\Big),
\end{eqnarray*}
and
\begin{eqnarray*}
~ &&[\Pi,\Pi]_{\rm II}(x,y,z,w,t)\\
~ &=&2\Big(\big(\Pi\circ\Pi\big)_{\rm II}(x,y,z,w)\Big)\\
~ &=&2\Big(\omega(x,y,\omega(z,w,t))-\omega(z,w,\omega(x,y,t))-\omega(\omega(x,y,z),w,t)-\omega(z,\omega(x,y,w),t)\Big),
\end{eqnarray*}
which proves the last statement. This completes the proof.
\end{proof}

\begin{rmk}
We have to point out that any Maurer-Cartan element $\Pi=(\pi,\omega)\in \frkC^1(\g,\g)$ in the graded Lie algebra $\frkC^\bullet(\g,\g)$ does not precisely correspond to a Lie-Yamaguti algebra structure on $\g$. By the proof of Theorem \ref{fund}, we know that $\Pi$ is a Maurer-Cartan element if and only if $\Pi$ satisfies Eqs. \eqref{LY3} and \eqref{fundamental}, i.e., for all $x,y,z,w,t\in \g$, $\Pi$ satisfies
\begin{eqnarray*}
\omega(x,y,\pi(z,w))&=&\pi(\omega(x,y,z),w)+\pi(z,\omega(x,y,w)),\\
~\omega(x,y,\omega(z,w,t))&=&\omega(\omega(x,y,z),w,t)+\omega(z,\omega(x,y,w),t)+\omega(z,w,\omega(x,y,t)).
\end{eqnarray*}
This reveals to us that case of Lie-Yamaguti algebras is more complicated than that of Lie algebras, Leibniz algebras or even $3$-Lie algebras. See \cite{GLST} for more details about Maurer-Cartan characterizations of other algebraic structures.
\end{rmk}

Let $\Pi$ define a \LYA ~structure on the vector space $\g$. It follows from the graded Jacobi identity that $d_\Pi:=[\Pi,\cdot]_{\rm LieY}$ is a differential on $(\frkC^\bullet(\g,\g),[\cdot,\cdot]_{\rm LieY})$. More precisely, we have
\begin{thm}\label{fund2}
Let $(\g,\pi,\omega)$ be a Lie-Yamaguti algebra. Then with the above notations, the triple $(\frkC^\bullet(\g,\g),[\cdot,\cdot]_{\rm LieY},d_{\Pi})$ is a differential graded Lie algebra, where $d_{\Pi}$ is defined to be
$$d_{\Pi}:=[\Pi,\cdot]_{\rm LieY}.$$
If moreover, for an element $\Pi'\in \frkC^1(\g,\g)$, $\Pi+\Pi'$ defines a Lie-Yamaguti algebra structure on $\g$, then $\Pi'$ is a Maurer-Cartan element of the differential graded Lie algebra $(\frkC^\bullet(\g,\g),[\cdot,\cdot]_{\rm LieY},d_{\Pi})$.
\end{thm}

From Proposition \ref{semi} and Theorem \ref{fund}, representations of Lie-Yamaguti algebras can be characterized as Maurer-Cartan elements in a graded Lie algebra.
More precisely, we have

\begin{pro}
With the above notations, if $(V;\rho,\mu)$ is a representation of Lie-Yamaguti algebra $(\g,\pi,\omega)$, then $\overline\Theta\in \frkC^1(\g\oplus V,\g \oplus V)$ is a Maurer-Cartan element of the differential graded Lie algebra $(\frkC^\bullet(\g\oplus V,\g\oplus V),[\cdot,\cdot]_{\rm LieY},d_{\overline\Pi})$.
\end{pro}

At the end of this section, we give the relationship between the coboundary operator of a Lie-Yamaguti algebra associated to the adjoint representation and the differential $d_{\Pi}$ induced by the Maurer-Cartan element.

\begin{thm}\label{diffe}
Let $(\g,\pi,\omega)$ be a Lie-Yamaguti algebra and $\delta:\frkC^n_{\mathsf{LieY}}(\g,\g)\to \frkC^{n+1}_{\mathsf{LieY}}(\g,\g)$ be the coboundary map associated to the adjoint representation. Then we have
\begin{eqnarray}
\label{cohomo1}\delta(F)&=&(-1)^{n}d_{\Pi}(F)=(-1)^{n}[\Pi,F]_{\rm LieY}, \quad \forall F\in \frkC^n(\g,\g) ~(n\geqslant 1),\\
\label{cohomo2}\delta(f)&=&d_{\Pi}(f)=[\Pi,f]_{\rm LieY},\quad \forall f\in \frkC^0(\g,\g)=\Hom(\g,\g).
\end{eqnarray}
\end{thm}
\begin{proof}
For all $F=(f,g)\in \frkC^n(\g,\g)~(n\geqslant 1)$ and for all $\frkX_1,\cdots,\frkX_{n+1}\in \wedge^2\g,~x\in \g$, we compute that
\begin{eqnarray*}
~ &&[\Pi,F]_{\rm I}(\frkX_1,\cdots,\frkX_{n+1})\\
~ &=&\Big(\big(\Pi\circ F\big)_{\rm I}-(-1)^n\big(F\circ\Pi\big)_{\rm I}\Big)(\frkX_1,\cdots,\frkX_{n+1})\\
~ &\stackrel{\eqref{gradedbra1}}{=}&\sum_{\sigma\in \mathbb S_{(1,n)}\atop \sigma(n+1)=n+1}(-1)^n sign(\sigma)\Courant{\frkX_{\sigma(1)},f(\frkX_{\sigma(2)},\cdots,\frkX_{\sigma(n+1)})}\\
~ &&+[x_{n+1},g(\frkX_1,\cdots,\frkX_n,y_{n+1})]-[y_{n+1},g(\frkX_1,\cdots,\frkX_n,x_{n+1})]\\
~ &&-(-1)^n\Big((-1)^ng(\frkX_{\sigma(1)},\cdots,\frkX_{\sigma(n)},[x_{n+1},y_{n+1}])\\
~ &&+\sum_{k=1}^n(-1)^{k-1}\sum_{\sigma\in \mathbb S_{(k-1,1)}}sign(\sigma)f(\frkX_{\sigma(1)},\cdots,\frkX_{\sigma(k-1)},
x_{k+1}\wedge\Courant{\frkX_{\sigma(k)},y_{k+1}},\frkX_{k+2},\cdots,\frkX_{n+1})\\
~ &&+\sum_{k=1}^n(-1)^{k-1}\sum_{\sigma\in \mathbb S_{(k-1,1)}}sign(\sigma)f(\frkX_{\sigma(1)},\cdots,\frkX_{\sigma(k-1)},
\Courant{\frkX_{\sigma(k)},x_{k+1}}\wedge y_{k+1},\frkX_{k+2},\cdots,\frkX_{n+1})\Big)\\
~ &=&(-1)^{n}\Big((-1)^n\big([x_{n+1},g(\frkX_1,\cdots,\frkX_n,y_{n+1})]-[y_{n+1},g(\frkX_1,\cdots,\frkX_n,x_{n+1})]\\
~ &&-g(\frkX_1,\cdots,\frkX_n,[x_{n+1},y_{n+1}])\big)\\
~ &&+\sum_{k=1}^{n}(-1)^{k+1}\Courant{\frkX_k,f(\frkX_1,\cdots,\hat{\frkX_k},\cdots,\frkX_{n+1})}\\
~ &&+\sum_{1\leqslant k<l\leqslant n+1}(-1)^{k}f(\frkX_1,\cdots,\hat{\frkX_k},\cdots,\frkX_k\circ\frkX_l,\cdots,\frkX_{n+1})\Big)\\
~ &\stackrel{\eqref{cohomolo1}}{=}&(-1)^{n}\delta_{\rm I}(F)(\frkX_1,\cdots,\frkX_{n+1}),
\end{eqnarray*}
and
\begin{eqnarray*}
~ &&[\Pi,F]_{\rm II}(\frkX_1,\cdots,\frkX_{n+1},x)\\
~ &=&\Big(\big(\Pi\circ F\big)_{\rm II}-(-1)^n\big(F\circ\Pi\big)_{\rm II}\Big)(\frkX_1,\cdots,\frkX_{n+1},x)\\
~ &\stackrel{\eqref{gradedbra3}}{=}&\sum_{\sigma\in \mathbb S_{(1,n)}}(-1)^n sign(\sigma)\Courant{\frkX_{\sigma(1)},g(\frkX_{\sigma(2)},\cdots,\frkX_{\sigma(n+1)},x)}\\
~ &&+\Courant{g(\frkX_1,\cdots,\frkX_n,x_{n+1}),y_{n+1},x}+\Courant{x_{n+1},g(\frkX_{1},\cdots,\frkX_n,y_{n+1}),x}\\
~ &&-(-1)^n\Big(\sum_{\sigma\in \mathbb S_{(n,1)}}(-1)^n sign(\sigma)g(\frkX_{\sigma(1)},\cdots,\frkX_{\sigma(n)},\Courant{\frkX_{\sigma(n+1)},x})\\
~ &&+\sum_{k=1}^n(-1)^{k-1}\sum_{\sigma\in \mathbb S_{(k-1,1)}}sign(\sigma)g(\frkX_{\sigma(1)},\cdots,\frkX_{\sigma(k-1)},\Courant{\frkX_{\sigma(k)},x_{k+1}}
\wedge y_{k+1},\frkX_{k+2},\cdots,\frkX_{n+1},x)\\
~ &&+\sum_{k=1}^n(-1)^{k-1}\sum_{\sigma\in \mathbb S_{(k-1,1)}}sign(\sigma)g(\frkX_{\sigma(1)},\cdots,\frkX_{\sigma(k-1)},x_{k+1}\wedge\Courant{\frkX_{\sigma{(k)}},y_{k+1}},\frkX_{k+2},\cdots,\frkX_{n+1},x)\Big)\\
~ &=&(-1)^{n}\Big((-1)^n\big(\Courant{g(\frkX_1,\cdots,\frkX_n,x_{n+1}),y_{n+1},x}-\Courant{g(\frkX_1,\cdots,\frkX_n,y_{n+1}),x_{n+1},x}\big)\\
~ &&+\sum_{k=1}^{n+1}(-1)^{k+1}\Courant{x_k,y_k,g(\frkX_1,\cdots,\hat{\frkX_k},\cdots,\frkX_{n+1},x)}\\
~ &&+\sum_{1\leqslant k<l\leqslant n+1}(-1)^kg(\frkX_1,\cdots,\hat{\frkX_k},\cdots,\frkX_k\circ\frkX_l,\cdots,\frkX_{n+1},x)\\
~ &&+\sum_{k=1}^{n+1}(-1)^kg(\frkX_1,\cdots,\hat{\frkX_k},\cdots,\frkX_{n+1},\Courant{x_k,y_k,x})\Big)\\
~ &\stackrel{\eqref{cohomolo2}}{=}&(-1)^{n}\delta_{\rm II}(F)(\frkX_1,\cdots,\frkX_{n+1},x).
\end{eqnarray*}
Thus, we have
$$[\Pi,F]_{\rm I}=(-1)^n\delta_{\rm I}(F)\quad \text{and}\quad [\Pi,F]_{\rm II}=(-1)^n\delta_{\rm II}(F),\quad \forall F\in \frkC^n(\g,\g)~(n\geqslant 1),$$
which proves \eqref{cohomo1}. Similarly, for any $f\in \Hom(\g,\g)$ and for all $x,y,z\in \g$, we have
\begin{eqnarray*}
~ &&[\Pi,f]_{\rm I}(x,y)\\
~ &=&\Big(\Pi\circ f\Big)_{\rm I}(x,y)-\Big(f\circ \Pi\Big)_{\rm I}(x,y)\\
~ &\stackrel{\eqref{gradedbra4},\eqref{gradedbra5}}{=}&[x,f(y)]+[f(x),y]-f([x,y])\\
~ &\stackrel{\eqref{1cochain}}{=}&\delta_{\rm I}(f)(x,y),
\end{eqnarray*}
and
\begin{eqnarray*}
~ &&[\Pi,f]_{\rm II}(x,y,z)\\
~ &\stackrel{\eqref{gradedbra6},\eqref{gradedbra2}}{=}&\Courant{f(x),y,z}+\Courant{x,f(y),z}+\Courant{x,y,f(z)}-f(\Courant{x,y,z})\\
~ &\stackrel{\eqref{2cochain}}{=}&\delta_{\rm II}(f)(x,y,z).
\end{eqnarray*}
Thus we have
$$[\Pi,f]_{\rm I}=\delta_{\rm I}(f)\quad\text{and}\quad [\Pi,f]_{\rm II}=\delta_{\rm II}(f), \quad\forall f\in \frkC^0(\g,\g)=\Hom(\g,\g),$$
which proves \eqref{cohomo2}. The conclusion thus follows.
\end{proof}

\section{Maurer-Cartan characterizations and Cohomology of relative Rota-Baxter operators}
In this section, we determine an $L_\infty$-algebra as controlling algebra for relative Rota-Baxter operators and then establish their cohomology theory. Finally, we reveal the relation between cohomology and differential in the controlling algebra.
\subsection{Maurer-Cartan characterizations for relative Rota-Baxter operators}
In this subsection, we use the Voronov's derived bracket to construct an $L_{\infty}$-algebra, which characterizes the relative Rota-Baxter operators on a Lie-Yamaguti algebra as Maurer-Cartan elements. At the beginning, we recall some notions and conclusions.
The notion of an $L_{\infty}$-algebra was introduced by Stasheff in \cite{Stasheff}. See also \cite{Lada2,Lada1} for more details.

\begin{defi}
An {\bf $L_{\infty}$-algebra} is a $\mathbb Z$-graded vector space $\g=\oplus_{k\in \mathbb Z}\g_k$ equipped with a collection of multilinear maps $l_k:\otimes^k\g \to \g~(k\geqslant 1)$ of degree $1$ with the property that, for all homogeneous elements $x_1,\cdots,x_n \in \g$, we have
\begin{itemize}
\item[\rm (i)] for all $\sigma\in S_n$,
$$l_n(x_{\sigma(1)},\cdots, x_{\sigma(n)})=\varepsilon(\sigma)l_n(x_1,\cdots,x_n).$$
\item[\rm (ii)] for all $n\geqslant 1$,
$$\sum_{i=1}^n\sum_{\sigma\in \mathbb S_{(i,n-1)}}\varepsilon(\sigma)l_{n-i+1}(l_i(x_{\sigma(1)},\cdots,x_{\sigma(i)}),x_{\sigma(i+1)},\cdots,x_{\sigma(n)})=0,$$
\end{itemize}
where $\varepsilon(\sigma)$ means the Koszul sign: $l_n(x_1,\cdots,x_i,x_{i+1},\cdots,x_n)=(-1)^{|x_i||x_{i+1}|}l_n(x_1,\cdots,x_{i+1},x_i,\cdots,x_n)$. In the sequel, we denote an $L_\infty$-algebra by $(\g,\{l_k\}_{k=1}^\infty)$.
\end{defi}

\begin{defi}
Let $(\g,\{l_k\}_{k=1}^\infty)$ be an $L_{\infty}$-algebra.
\begin{itemize}
\item[\rm (i)] A degree $0$ element $\alpha \in \g_0$ is called a {\bf Maurer-Cartan element} of the $L_{\infty}$-algebra $(\g,\{l_k\}_{k=1}^\infty)$ if $\alpha$ satisfies  the following Maurer-Cartan equation:
\begin{eqnarray*}
\sum_{k=1}^\infty \frac{1}{k!}l_k(\underbrace{\alpha,\cdots,\alpha}_k)=0.
\end{eqnarray*}
\item[\rm (ii)] A Maurer-Cartan element $\alpha\in \g_0$ of $(\g,\{l_k\}_{k=1}^\infty)$ is called  {\bf strict} if
$$l_k(\alpha,\cdots,\alpha)=0,\quad k\geqslant 1.$$
\end{itemize}
\end{defi}

Let us recall the method by which we use the Voronov's derived brackets to construct an $L_{\infty}$-algebra.

\begin{defi}(\cite{Voronov})
A {\bf $V$-data} is a quadruple $(L,\h,P,\Delta)$ consisting of the following items:
\begin{itemize}
\item $(L,[\cdot,\cdot])$ is a graded Lie algebra;
\item $\h$ is an abelian graded Lie subalgebra of $(L,[\cdot,\cdot])$;
\item $P:L \to L$ is a projection, i.e. $P\circ P=P$, whose image is $\h$ and its kernel is a graded Lie subalgebra of $(L,[\cdot,\cdot])$;
\item $\Delta$ is an element in $\ker(P)^1$ such that $[\Delta,\Delta]=0$.
\end{itemize}
\end{defi}

\begin{thm}{\rm (\cite{Voronov})}\label{construct}
Let $(L,\h,P,\Delta)$ be a $V$-data. Define multilinear maps $l_k:\otimes^k\h\to \h$ to be
\begin{eqnarray}
l_k(x_1,\cdots,x_k)=P[\cdots[[\Delta,x_1],x_2],\cdots,x_k],\label{derived}
\end{eqnarray}
for all homogeneous elements $x_1,\cdots,x_k\in \h$. Then $(\h,\{l_k\}_{k=1}^\infty)$  is an $L_\infty$-algebra, where $\{l_k\}_{k=1}^\infty$ is called the {\bf higher derived brackets} of the $V$-data $(L,\h,P,\Delta)$.
\end{thm}

\begin{defi}{\rm\cite{SZ1}}
Let $(\g,[\cdot,\cdot],\Courant{\cdot,\cdot,\cdot})$ be a Lie-Yamaguti algebra and $(V;\rho,\mu)$ a representation of $\g$. A linear map $T:V\to \g$ is called a {\bf relative Rota-Baxter operator} on $\g$ with respect to $(V;\rho,\mu)$ if $T$ satisfies
\begin{eqnarray}
~\label{Ooperator1}[Tu,Tv]&=&T\Big(\rho(Tu)v-\rho(Tv)u\Big),\\
~\label{Ooperator2}\Courant{Tu,Tv,Tw}&=&T\Big(D_{\rho,\mu}(Tu,Tv)w+\mu(Tv,Tw)u-\mu(Tu,Tw)v\Big), \quad \forall u,v,w \in V.
\end{eqnarray}
\end{defi}
\emptycomment{
\begin{rmk}
If a Lie-Yamaguti algebra $(\g,[\cdot,\cdot],\Courant{\cdot,\cdot,\cdot})$ with a representation $(V;\rho,\mu)$ reduces to a Lie triple system $(\g,\Courant{\cdot,\cdot,\cdot})$ with a  representation   $(V;\mu)$, we   obtain the notion of  a {\bf relative Rota-Baxter operator on a Lie triple system}, i.e. the following equation holds:
\begin{eqnarray*}
\Courant{Tu,Tv,Tw}=T\Big(D_\mu(Tu,Tv)w+\mu(Tv,Tw)u-\mu(Tu,Tw)v\Big), \quad \forall u,v,w \in V,
\end{eqnarray*}
where $D_\mu(x,y):=\mu(y,x)-\mu(x,y)$ for any $x,y \in \g.$ Thus, all the results given in the sequel can be adapted to the Lie triple system context.
\end{rmk}
}

In the sequel, let $\g$ and $V$ be vector spaces, and we denote the elements in $\g$ by $x,y,z$, and we denote the elements in $V$ by $u,v,w$. Consider the graded vector subspace $C^*(V,\g)=\oplus_{p\geqslant 0}C^p(V,\g)$ of $\frkC^\bullet(\g\oplus V,\g\oplus V)$, where
\begin{eqnarray*}
C^p(V,\g)=
\begin{cases}
\Hom(\underbrace{\wedge^2V\otimes\cdots\otimes\wedge^2V}_p,\g)\oplus\Hom(\underbrace{\wedge^2V\otimes\cdots\otimes\wedge^2V}_p\otimes V,\g),& p\geqslant 1,\\
\Hom(V,\g),&p=0.
\end{cases}
\end{eqnarray*}

Recall that $\Pi+\Theta=(\pi+\rho,\omega+\mu)\in \frkC^1(\g\oplus V,\g\oplus V)$. For all $P\in C^p(V,\g),~Q\in C^q(V,\g)$, and $R\in C^r(V,\g)~(p,q,r\geqslant 0)$, we define a bilinear operation
$$l_2:C^p(V,\g)\times C^q(V,\g) \to C^{p+q+1}(V,\g)$$
and a trilinear operation
$$l_3:C^p(V,\g)\times C^q(V,\g) \times C^r(V,\g) \to C^{p+q+r+1}(V,\g)$$
to be
\begin{eqnarray}
\label{multi1}l_2(P,Q)&=&[[\Pi+\Theta,P]_{\mathsf{LieY}},Q]_{\mathsf{LieY}},
\end{eqnarray}
and
\begin{eqnarray}
\label{multi2}l_3(P,Q,R)&=&[[[\Pi+\Theta,P]_{\mathsf{LieY}},Q]_{\mathsf{LieY}},R]_{\mathsf{LieY}},
\end{eqnarray}
respectively. Note that the standard symbols for $l_2$ and $l_3$ are $\Big((l_2)_{\rm I},(l_2)_{\rm II}\Big)$ and $\Big((l_3)_{\rm I},(l_3)_{\rm II}\Big)$ respectively, if we emphasize their components.

Now we give the main theorem in this section which gives the Maurer-Cartan characterization for relative Rota-Baxter operators on Lie-Yamaguti algebras.

\begin{thm}\label{main}
With the above notations, $(C^*(V,\g),l_2,l_3)$ is an $L_\infty$-algebra. Moreover, its strict Maurer-Cartan elements are precisely the relative Rota-Baxter operators on Lie-Yamaguti algebra $(\g,\pi,\omega)$ with respect to the representation $(V;\rho,\mu)$.
\end{thm}
\begin{proof}
Let $(V;\rho,\mu)$ be a representation  of the Lie-Yamaguti algebra $(\g,\pi,\omega)$. Then we have
\begin{itemize}
\item a graded Lie algebra $(\frkC^\bullet(\g\oplus V,\g\oplus V),[\cdot,\cdot]_{\mathsf{LieY}})$;
\item an abelian subalgebra $C^*(V,\g)$;
\item $\pr_\g$ is a projection onto $C^*(V,\g)$;
\item $\Delta=\Pi+\Theta$.
\end{itemize}
These items listed above forms a $V$-data. By Theorem \ref{construct}, $(C^*(V,\g),\{l_k\}_{k=1}^\infty)$ is an $L_\infty$-algebra, where multilinear maps ${l_k}'s$ are given by \eqref{derived}. But we note that
\begin{eqnarray*}
~\pr_\g[\Pi+\Theta,P]_{\mathsf{LieY}}=0,\\
~\pr_\g[[\Pi+\Theta,P]_{\mathsf{LieY}},Q]_{\mathsf{LieY}}\in C^*(V,\g),\\
~\pr_\g[[[\Pi+\Theta,P]_{\mathsf{LieY}},Q]_{\mathsf{LieY}},R]_{\mathsf{LieY}}\in C^*(V,\g),
\end{eqnarray*}
for all $P\in C^p(V,\g),~Q\in C^q(V,\g)$ and $R\in C^r(V,\g)$. Thus, the nonzero multilinear maps of $L_\infty$-algebra $(C^*(V,\g),\{l_k\}_{k=1}^\infty)$ are only $l_2$ and $l_3$.

For a degree $0$ element $T\in \Hom(V,\g)$, by \eqref{multi1} and \eqref{multi2}, we have for all $u,v,w \in V$,
\begin{eqnarray*}
~ &&(l_2)_{\rm I}(T,T)(u,v)\\
~ &=&[[\Pi+\Theta,T]_{\mathsf{LieY}},T]_{\rm I}(u,v)\\
~ &=&\Big([\Pi+\Theta,T]_{\mathsf{LieY}}\circ T\Big)_{\rm I}(u,v)-\Big(T\circ[\Pi+\Theta,T]_{\mathsf{LieY}} \Big)_{\rm I}(u,v)\\
~ &=&[\Pi+\Theta,T]_{\rm I}(Tu,v)+[\Pi+\Theta,T]_{\rm I}(u,Tv)\\
~ &=&(\Pi+\Theta)(Tu,Tv)-\Big(T\circ(\Pi+\Theta) \Big)_{\rm I}(Tu,v)\\
~ &&+(\Pi+\Theta)(Tu,Tv)-\Big(T\circ(\Pi+\Theta) \Big)_{\rm I}(u,Tv)\\
~ &=&2\Big([Tu,Tv]-T\big(\rho(Tu)v-\rho(Tv)u\big)\Big),
\end{eqnarray*}
and
\begin{eqnarray*}
~ &&(l_3)_{\rm II}(T,T,T)(u,v,w)\\
~ &=&[[[\Pi+\Theta,T]_{\mathsf{LieY}},T]_{\mathsf{LieY}},T]_{\rm II}(u,v,w)\\
~ &=&[[\Pi+\Theta,T]_{\mathsf{LieY}},T]_{\rm II}(Tu,v,w)+[[\Pi+\Theta,T]_{\mathsf{LieY}},T]_{\rm II}(u,Tv,w)\\
~ &&+[[\Pi+\Theta,T]_{\mathsf{LieY}},T]_{\rm II}(u,v,Tw)-\Big(T\circ [[\Pi+\Theta,T]_{\mathsf{LieY}},T]_{\mathsf{LieY}}\Big)_{\rm II}(u,v,w)\\
~ &=&[\Pi+\Theta,T]_{\rm II}(Tu,Tv,w)+[\Pi+\Theta,T]_{\rm II}(Tu,v,Tw)\\
~ &&-\Big(T\circ [\Pi+\Theta,T]_{\mathsf{LieY}}\Big)_{\rm II}(Tu,v,w)+[\Pi+\Theta,T]_{\rm II}(Tu,Tv,w)\\
~ &&+[\Pi+\Theta,T]_{\rm II}(u,Tv,Tw)-\Big(T\circ [\Pi+\Theta,T]_{\mathsf{LieY}}\Big)_{\rm II}(u,Tv,w)\\
~ &&+[\Pi+\Theta,T]_{\rm II}(Tu,v,Tw)+[\Pi+\Theta,T]_{\rm II}(u,Tv,Tw)\\
~ &&-\Big(T\circ [\Pi+\Theta,T]_{\mathsf{LieY}}\Big)_{\rm II}(u,v,Tw)-\Big(T\circ [\Pi+\Theta,T]_{\mathsf{LieY}}\Big)_{\rm II}(Tu,v,w)\\
~ &&-\Big(T\circ [\Pi+\Theta,T]_{\mathsf{LieY}}\Big)_{\rm II}(u,Tv,w)-\Big(T\circ [\Pi+\Theta,T]_{\mathsf{LieY}}\Big)_{\rm II}(u,v,Tw)\\
~ &=&6\Big(\Courant{Tu,Tv,Tw}-T\big(D(Tu,Tv)w+\mu(Tv,Tw)u-\mu(Tu,Tw)v\big)\Big).
\end{eqnarray*}
To sum up, we compute the nonzero multilinear maps $l_2$  and $l_3$ as follows:
\begin{eqnarray*}
\begin{cases}
(l_2)_{\rm I}(T,T)(u,v)=2\Big([Tu,Tv]-T\big(\rho(Tu)v-\rho(Tv)u\big)\Big),\\
(l_2)_{\rm II}(T,T)(u,v,w)=0,\\
(l_3)_{\rm I}(T,T,T)(u,v)=0,\\
(l_3)_{\rm II}(T,T,T)(u,v,w)=6\Big(\Courant{Tu,Tv,Tw}-T\big(D(Tu,Tv)w+\mu(Tv,Tw)u-\mu(Tu,Tw)v\big)\Big).
\end{cases}
\end{eqnarray*}
Thus we have that a linear map $T:V\to \g$ is a relative Rota-Bxter operator on the Lie-Yamaguti algebra $(\g,\pi,\omega)$ with respect to the representation $(V;\rho,\mu)$ if and only if $T$ is a strict Maurer-Cartan element of the $L_\infty$-algebra $(C^*(V,\g),l_2,l_3)$. This completes the proof.
\end{proof}

\begin{rmk}
From Theorem \ref{main}, we notice that for a relative Rota-Baxter operator on a \LYA, on the one hand, its realization of Maurer-Cartan characterization is an {\em $L_\infty$-algebra} owning two nonzero multilinear maps: $l_2$ and $l_3$. In \cite{TBGS,THS}, authors showed the realizations of Maurer-Cartan characterization for relative Rota-Baxter operators on Lie algebras and $3$-Lie algebras are graded Lie algebras and Lie $3$-algebras respectively, both of which can be seen as a special $L_\infty$-algebra (whose nonzero multilinear maps are $l_2$ and $l_3$ respectively). Whichever is a graded Lie algebra or a Lie $3$-algebra, it (treated as $L_\infty$-algebra) owns one nonzero bracket. On the other hand, a relative Rota-Baxter operator on a \LYA~ corresponds to a {\em strict} Maurer-Cartan element in its $L_\infty$-algebra, not a Maurer-Cartan element. The above two points distinguish from the case of relative Rota-Baxter operators on Lie algebras and on $3$-Lie algebras \cite{TBGS,THS,T.S2}, and from proof of Theorem \ref{main}, it reveals the fact that the case of \LYA s is more complicated than other algebras.
\end{rmk}

Let $T$ be a relative Rota-Baxter operator on a Lie-Yamaguti algebra $(\g,\pi,\omega)$ with respect to a representation $(V;\rho,\mu)$. Since $T\in C^0(V,\g)=\Hom(V,\g)$ is a (strict) Maurer-Cartan element of the $L_\infty$-algebra $(C^*(V,\g),l_2,l_3)$, by \cite{Get}, we have a twisted $L_\infty$-algebra $(C^*(V,\g),l_1^T,l_2^T,l_3^T)$ as follows:
\begin{eqnarray}
~l_1^T(P)&=&l_2(T,P)+\half l_3(T,T,P),\label{diffs}\\
~l_2^T(P,Q)&=&l_2(P,Q)+l_3(T,P,Q),\label{bracket}\\
~l_3^T(P,Q,R)&=&l_3(P,Q,R),\label{homotopy}\\
~l_k^T&=&0,\quad k\geqslant 4,
\end{eqnarray}
for all $P\in C^p(V,\g),~Q\in C^q(V,\g)$ and $R\in C^r(V,\g)$.

\begin{thm}\label{deformation}
Let $T:V\to \g$ be a relative Rota-Baxter operator on a Lie-Yamaguti algebra $(\g,[\cdot,\cdot],\Courant{\cdot,\cdot,\cdot})$ with respect to a representation $(V;\rho,\mu)$. Then for a linear map $T':V\to \g$, $T+T'$ is still a relative Rota-Baxter operator on $\g$ if and only if $T'$ is a Maurer-Cartan element of the twisted $L_\infty$-algebra $(C^*(V,\g),l_1^T,l_2^T,l_3^T)$, i.e., $T'$ satisfies the following Maurer-Cartan equation:
$$l_1^T(T')+\half l_2^T(T',T')+\frac{1}{3!}l_3^T(T',T',T')=0.$$
\end{thm}
\begin{proof}
By Theorem \ref{main}, $T+T'$ is a relative Rota-Baxter operator if and only if it is a strict Maurer-Cartan element, i.e.,
\begin{eqnarray}\label{Maurer}
\begin{cases}
l_2(T+T',T+T')=0,\\
l_3(T+T',T+T',T+T')=0.
\end{cases}
\end{eqnarray}
Since $T$ is a relative Rota-Baxter operator on $(\g,[\cdot,\cdot],\Courant{\cdot,\cdot,\cdot})$, then we have
\begin{eqnarray*}
\begin{cases}
l_2(T,T)=0,\\
l_3(T,T,T)=0.
\end{cases}
\end{eqnarray*}
Since the degrees of both $T$ and $T'$ are $0$, expanding Eqs.\eqref{Maurer} yields that
\begin{eqnarray}
\begin{cases}
2l_2(T,T')+l_2(T',T')=0,\\
3l_3(T,T,T')+3l_3(T,T',T')+l_3(T',T',T')=0.\label{MC}
\end{cases}
\end{eqnarray}
Then we obtain that
\begin{eqnarray*}
~ &&l_1^T(T')+\half l_2^T(T',T')+\frac{1}{3!}l_3^T(T',T',T')\\
~ &\stackrel{\eqref{diffs},\eqref{bracket},\eqref{homotopy}}{=}&l_2(T,T')+\half l_3(T,T,T')+\half\Big(l_2(T',T')+l_3(T,T',T')\Big)+\frac{1}{3!}l_3(T',T',T')\\
~ &=&\Big(l_2(T,T')+\half l_2(T',T')\Big)+\half\Big(l_3(T,T,T')+l_3(T,T',T')+\frac{1}{3}l_3(T',T',T')\Big)\\
~ &\stackrel{\eqref{MC}}{=}&0,
\end{eqnarray*}
which implies that $T'$ is a Maurer-Cartan element of the twisted $L_\infty$-algebra $(C^*(V,\g),l_1^T,l_2^T,l_3^T)$. This completes the proof.
\end{proof}

At the end of this subsection, we give two examples of Rota-Baxter operators on \LYA s.

\begin{ex}
Let $\g=C[0,1]$ endowed with the following operations
\begin{eqnarray*}
[f,g](x)&=&f(x)g(x)-g(x)f(x),\\
\Courant{f,g,h}(x)&=&f(x)g(x)h(x)-g(x)f(x)h(x)-h(x)f(x)g(x)+h(x)g(x)f(x),\quad \forall x\in [0,1],
\end{eqnarray*}
for all $ f,g,h\in \g$. Then $(\g,[\cdot,\cdot],\Courant{\cdot,\cdot,\cdot})$ forms a \LYA. For all $f\in C[0,1]$, the {\em integral operator} $R:\g\longrightarrow\g$ defined to be
$$R(f)(x):=\int_0^xf(t){\mathrm d}t,\quad\forall x\in [0,1]$$
is a Rota-Baxter operator on $\g$. For the integral operator $R$, define $l_2$ and $l_3$ to be
$$
\begin{cases}
(l_2)_{\rm I}(R,R)(f,g)=2\Big([R(f),R(g)]-R([R(f),g]-[f,R(g)])\Big),\\
(l_2)_{\rm II}(R,R)(f,g,h)=0,\\
(l_3)_{\rm I}(R,R,R)(f,g)=0,\\
(l_3)_{\rm II}(R,R,R)(f,g,h)=6\Big(\Courant{R(f),R(g),R(h)}-R\big(\Courant{R(f),R(g),h}-\Courant{f,R(g),R(h)}-\Courant{R(f),g,R(h)}\big)\Big),
\end{cases}
$$
for all $f,g,h\in \g$. Obviously, $R$ is a strict Maurer-Cartan element of the $L_\infty$-algebra $(C^*(\g,\g),l_2,l_3)$ by Theorem \ref{main}, where $l_2$ and $l_3$ are defined by \eqref{multi1} and \eqref{multi2}. Moreover, for another linear operator $R':\g\longrightarrow\g$, $R+R'$ is still a Rota-Baxter on $\g$ if and only if $R'$ satisfies the following Maurer-Cartan equation:
$$l_1^R(R')+\frac{1}{2}l_2^R(R',R')+\frac{1}{3!}l_3^R(R',R',R')=0,$$
where $l_1^R,~l_2^R$ and $l_3^R$ are defined by Eqs. \eqref{diffs}-\eqref{homotopy} respectively.
\end{ex}

\begin{ex}
Let $(A,\cdot)$ be an associative algebra and $A[[\nu]]$ an algebra of formal series with coefficients in $A$. Define operations $[\cdot,\cdot]$ and $\Courant{\cdot,\cdot,\cdot}$ to be
\begin{eqnarray*}
[a_i\nu^i,a_j\nu^j]&:=&[a_i,a_j]\nu^{i+j},\\
~\Courant{a_i\nu^i,a_j\nu^j,a_k\nu^k}&:=&[[a_i,a_j],a_k]\nu^{i+j+k},
\end{eqnarray*}
where $a_i,a_j,a_k\in A$ and $[\cdot,\cdot]$ on the right hand side means the commutator: $[x,y]=x\cdot y-y\cdot x$ for all $x,y\in A$, then $(A[[\nu]],[\cdot,\cdot],\Courant{\cdot,\cdot,\cdot})$ forms a \LYA. Define a formal integral operator $\Omega:A[[\nu]]\longrightarrow A[[\nu]]$ to be
$$\Omega(a_i\nu^i):=\int a_i\nu^id\nu=\frac{1}{i+1}a_i\nu^{i+1}, \quad \forall a_i\in A,$$
then $\Omega$ is a Rota-Baxter operator on $A[[\nu]]$. Then we define $l_2$ and $l_3$ by Eqs. \eqref{multi1} and \eqref{multi2} respectively to achieve an $L_\infty$-algebra $(C^*(A[[\nu]],A[[\nu]]),l_2,l_3)$, and consequently we obtain that $\Omega$ is a strict Maurer-Cartan element, i.e.,
$$l_2(\Omega,\Omega)=0, \quad\text{and}\quad l_3(\Omega,\Omega,\Omega)=0. $$
\end{ex}

\subsection{Cohomology of relative Rota-Baxter operators}
In this subsection, we establish the cohomology of relative Rota-Baxter operators on Lie-Yamaguti algebras, and then we show the relation between the differential of the twisting $L_\infty$-algebra and the coboundary operator. We write $D_{\rho,\mu}$ to emphasize that it is relative with $\rho$ and $\mu$ in this subsection. 

\begin{pro}{\cite{SZ1}}
Let $T:V\longrightarrow\g$ be a relative Rota-Baxter operator on a Lie-Yamaguti algebra $(\g,[\cdot,\cdot],\Courant{\cdot,\cdot,\cdot})$ with respect to $(V;\rho,\mu)$. Define
\begin{eqnarray*}
[u,v]_T&=&\rho(Tu)v-\rho(Tv)u,\\
\Courant{u,v,w}_T&=&D_{\rho,\mu}(Tu,Tv)w+\mu(Tv,Tw)u-\mu(Tu,Tw)v,\quad \forall u,v,w\in V.
\end{eqnarray*}
Then $(V,[\cdot,\cdot]_T,\Courant{\cdot,\cdot,\cdot}_T)$ is a Lie-Yamaguti algebra, which is called the {\bf sub-adjacent Lie-Yamaguti algebra}.
Thus $T$ is a Lie-Yamaguti algebra homomorphism\footnote{A Lie-Yamaguti homomorphism is a linear map between two \LYA s that preserves the \LYA ~structure.} from $(V,[\cdot,\cdot]_T,\Courant{\cdot,\cdot,\cdot}_T)$ to $(\g,[\cdot,\cdot],\Courant{\cdot,\cdot,\cdot})$.
\end{pro}

Next, we give a representation of the sub-adjacent \LYA ~$(V,[\cdot,\cdot]_T,\Courant{\cdot,\cdot,\cdot}_T)$ on the vector space $\g$.
Define two linear maps $\varrho:V\to \gl(\g)$ and $\varpi:\otimes^2V\to \gl(\g)$ to be
\begin{eqnarray}
\label{repre1}\varrho(u)x&:=&[Tu,x]+T\big(\rho(x)u\big),\\
\label{repre2}\varpi(u,v)x&:=&\Courant{x,Tu,Tv}-T\big(D_{\rho,\mu}(x,Tu)v-\mu(x,Tv)u\big), \quad \forall x\in \g,~u,v \in V.
\end{eqnarray}

Consequently, we give the precise formula of $D_{\varrho,\varpi}$ first.
\begin{pro}
Let $T$ be a relative Rota-Baxter operator on a Lie-Yamaguti algebra $(\g,[\cdot,\cdot],\Courant{\cdot,\cdot,\cdot})$ with respect to $(V;\rho,\mu)$. Then with the notations above, we have
\begin{eqnarray}
\label{repre3}D_{\varrho,\varpi}(u,v)x=\Courant{Tu,Tv,x}-T\Big(\mu(Tv,x)u-\mu(Tu,x)v\Big), \quad \forall u,v \in V, ~x\in \g.
\end{eqnarray}
\end{pro}
\begin{proof}
Since $T$ is a relative Rota-Baxter operator, for all $u,v\in V$ and $x\in \g$, by a direct computation, we have
\begin{eqnarray*}
~ &&D_{\varrho,\varpi}(u,v)x\\
&\stackrel{\eqref{rep}}{=}&\varpi(v,u)x-\varpi(u,v)x+[\varrho(u),\varrho(v)]x-\varrho([u,v]_T)x\\
~ &\stackrel{\eqref{repre1},\eqref{repre2}}{=}&\Courant{x,Tv,Tu}-T\Big(D_{\rho,\mu}(x,Tv)u-\mu(x,Tu)v\Big)-\Courant{x,Tu,Tv}+T\Big(D_{\rho,\mu}(x,Tu)v-\mu(x,Tv)u\Big)\\
~ &&+[Tu,[Tv,x]]+[Tu,T(\rho(x)v)]+T(\rho([Tv,x])u)+T(\rho(T(\rho(x)v)u))\\
~ &&-[Tv,[Tu,x]]-[Tv,T(\rho(x)u)]-T(\rho([Tu,x])v)-T(\rho(T(\rho(x)u)v))\\
~ &&-[T(\rho(Tu)v-\rho(Tv)u),x]-T(\rho(x)\rho(Tu)v)+T(\rho(x)\rho(Tv)u)\\
~ &\stackrel{\eqref{Ooperator1}}{=}&\Courant{x,Tv,Tu}-\Courant{x,Tu,Tv}+[Tu,[Tv,x]]-[Tv,[Tu,x]]-[[Tu,Tv],x]\\
~ &&-T\Big(D_{\rho,\mu}(x,Tv)u-\mu(x,Tu)v\Big)+T\Big(D_{\rho,\mu}(x,Tu)v-\mu(x,Tv)u\Big)\\
~ &&+T(\rho(Tu)\rho(x)v-\rho(T(\rho(x)v)u)-T(\rho(Tv)\rho(x)u-\rho(T(\rho(x)u)v)\\
~ &&+T(\rho([Tv,x])u)+T(\rho(T(\rho(x)v)u))-T(\rho([Tu,x])v)-T(\rho(T(\rho(x)u)v))\\
~ &&-T(\rho(x)\rho(Tu)v)+T(\rho(x)\rho(Tv)u)\\
~ &\stackrel{\eqref{LY1},\eqref{rep}}{=}&\Courant{Tu,Tv,x}-T\Big(\mu(Tv,x)u-\mu(Tu,x)v\Big).
\end{eqnarray*}
The conclusion thus follows.
\end{proof}

\begin{thm}\label{represent}
With the above notations, then $(\g;\varrho,\varpi)$ is a representation of the sub-adjacent Lie-Yamaguti algebra $(V,[\cdot,\cdot]_T,\Courant{\cdot,\cdot,\cdot}_T)$, where linear maps $\varrho,~\varpi$, and $D_{\varrho,\varpi}$ are given by Eqs. \eqref{repre1}-\eqref{repre3} respectively.
\end{thm}

We may prove Theorem \ref{represent} by checking that linear maps $\varrho$, $\varpi$, and $D_{\varrho,\varpi}$ satisfy conditions in Definition \ref{defi:representation}, but here we prove this proposition by another way. In order to do this, we should go back some notions in \cite{Sheng Zhao}. Recall that a Nijenhuis operator on a Lie-Yamaguti algebra $(\g,[\cdot,\cdot],\Courant{\cdot,\cdot,\cdot})$ is a linear map $N:\g\to \g$ satisfying
\begin{eqnarray*}
[Nx,Ny]&=&N\big([Nx,y]+[x,Ny]-N[x,y]\big),\\
~\Courant{Nx,Ny,Nz}&=&N\Big(\Courant{Nx,Ny,z}+\Courant{Nx,y,Nz}+\Courant{x,Ny,Nz}\\
~ &&-N\Courant{Nx,y,z}-N\Courant{x,Ny,z}-N\Courant{x,y,Nz}+N^2\Courant{x,y,z}\Big), \quad \forall x,y,z \in \g.
\end{eqnarray*}
Then we get a pair of deformed brackets $([\cdot,\cdot]_N,\Courant{\cdot,\cdot,\cdot}_N)$:
\begin{eqnarray}
\label{deform1}[x,y]_N&=&[Nx,y]+[x,Ny]-N[x,y],\\
\label{deform2}~\Courant{x,y,z}_N&=&\Courant{Nx,Ny,z}+\Courant{Nx,y,Nz}+\Courant{x,Ny,Nz}\\
\nonumber~ &&-N\Courant{Nx,y,z}-N\Courant{x,Ny,z}-N\Courant{x,y,Nz}+N^2\Courant{x,y,z}, \quad \forall x,y,z \in \g.
\end{eqnarray}
In \cite{Sheng Zhao}, authors showed that $(\g,[\cdot,\cdot]_N,\Courant{\cdot,\cdot,\cdot}_N)$ is a Lie-Yamaguti algebra and thus $N$ is a Lie-Yamaguti homomorphism from $(\g,[\cdot,\cdot]_N,\Courant{\cdot,\cdot,\cdot}_N)$ to $(\g,[\cdot,\cdot],\Courant{\cdot,\cdot,\cdot})$.

\vspace{5mm}
\emph{Proof of Theorem {\rm \ref{represent}:}}
It is direct to see that if $T:V \to \g$ is a relative Rota-Baxter operator on a Lie-Yamaguti algebra $(\g,[\cdot,\cdot],\Courant{\cdot,\cdot,\cdot})$ with respect to a representation $(V;\rho,\mu)$, then
$N_T=\begin{pmatrix}
0 & T\\
0 & 0
\end{pmatrix}$
is a Nijenhuis operator on the semidirect product Lie-Yamaguti algebra $\g\ltimes_{\rho,\mu} V$. Then by \eqref{deform1} and \eqref{deform2}, we deduce that there is a Lie-Yamaguti algebra structure on $V\oplus\g\cong\g \oplus V$ given by for all $x,y,z \in \g,~u,v,w\in V,$
\begin{eqnarray*}
~ &&[x+u,y+v]_{N_T}\\
~ &=&[N_T(x+u),y+v]_{\rho,\mu}+[x+u,N_T(y+v)]_{\rho,\mu}-N_T[x+u,y+v]_{\rho,\mu}\\
~ &=&[Tu,y+v]_{\rho,\mu}+[x+u,Tv]_{\rho,\mu}-N_T([x,y]+\rho(x)v-\rho(y)u)\\
~ &=&[Tu,y]+\rho(Tu)v+[x,Tv]-\rho(Tv)u-T(\rho(x)v-\rho(y)u)\\
~ &=&[u,v]_T+\varrho(u)y-\varrho(v)x,\\
~ &&\\
~&&\Courant{x+u,y+v,z+w}_{N_T}\\
~ &=&\Courant{N_T(x+u),N_T(y+v),z+w}_{\rho,\mu}+\Courant{N_{T}(x+u),y+v,N_{T}(z+w)}_{\rho,\mu}+\Courant{x+u,N_{T}(y+v),N_{T}(z+w)}_{\rho,\mu}\\
~ &&-N_T(\Courant{N_{T}(x+u),y+v,z+w}_{\rho,\mu}+\Courant{x+u,N_{T}(y+v),z+w}_{\rho,\mu})+\Courant{x+u,y+v,N_{T}(z+w)}_{\rho,\mu})\\
~ &=&\Courant{Tu,Tv,z+w}_{\rho,\mu}+\Courant{Tu,y+v,Tw}_{\rho,\mu}+\Courant{x+u,Tv,Tw}_{\rho,\mu}\\
~ &&-N_T(\Courant{Tu,y+v,z+w}_{\rho,\mu}+\Courant{x+u,Tv,z+w}_{\rho,\mu}+\Courant{x+u,y+v,Tw}_{\rho,\mu})\\
~ &=&\Courant{Tu,Tv,z}+D_{\rho,\mu}(Tu,Tv)w+\Courant{Tu,y,Tw}-\mu(Tu,Tw)v+\Courant{x,Tv,Tw}+\mu(Tv,Tw)u\\
~ &&-T\Big(D_{\rho,\mu}(Tu,y)w-\mu(Tu,z)v+D_{\rho,\mu}(x,Tv)w+\mu(Tv,z)u+\mu(y,Tw)u-\mu(x,Tw)v\Big)\\
~ &=&\Courant{u,v,w}_T+D_{\varrho,\varpi}(u,v)z+\varpi(v,w)x-\varpi(u,w)y,
\end{eqnarray*}
which implies that $(\g;\varrho,\varpi)$ is a representation of Lie-Yamaguti algebra $(V,[\cdot,\cdot]_T,\Courant{\cdot,\cdot,\cdot}_T)$. This finishes the proof.\qed

Next, we construct the $0$-cochians and the corresponding coboundary map.
For all $\frkX\in \wedge^2\g$, define $\delta(\frkX):V\to \g$ to be
\begin{eqnarray}
\label{delta}\delta(\frkX)v:=TD(\frkX)v-\Courant{\frkX,Tv}, \quad \forall v\in V.
\end{eqnarray}

\begin{pro}\label{0cocy}
Let $T$ be a relative Rota-Baxter operator on a Lie-Yamaguti algebra $(\g,[\cdot,\cdot,],\Courant{\cdot,\cdot,\cdot})$ with respect to a representation $(V;\rho,\mu)$. Then $\delta(\frkX)$ is a $1$-cocycle on the Lie-Yamaguti algebra $(V,[\cdot,\cdot]_T,\Courant{\cdot,\cdot,\cdot}_T)$ with coefficients in the representation $(\g;\varrho,\varpi)$.
\end{pro}
\begin{proof}
It is sufficient to show that both $\delta_{\rm I}^T(\delta(\frkX))$ and $\delta_{\rm II}^T(\delta(\frkX))$ all vanish. Indeed, for any $u,v,w\in V$, we have
\emptycomment{
{\footnotesize
\begin{eqnarray*}
~ &&\delta^T_{\rm I}\Big(\delta(\frkX)\Big)(u,v)\\
~ &\stackrel{\eqref{1cochain}}{=}&\varrho(u)\delta(\frkX)(v)-\varrho(v)\delta(\frkX)(u)-\delta(\frkX)([u,v]_T)\\
~ &\stackrel{\eqref{repre1}}{=}&[Tu,\delta(\frkX)(v)]+T(\rho(\delta(\frkX)(v))u)-[Tv,\delta(\frkX)(u)]-T(\rho(\delta(\frkX)(u))v)\\
~ &&-T(D_{\rho,\mu}(\frkX)[u,v]_T)+\Courant{\frkX,T[u,v]_T}\\
~ &\stackrel{\eqref{delta}}{=}&[Tu,TD_{\rho,\mu}(\frkX)v]-[Tu,\Courant{\frkX,Tv}]+T(\rho(T(D_{\rho,\mu}(\frkX)v))u)-T(\rho(\Courant{\frkX,Tv})u)\\
~ &&-[Tv,TD_{\rho,\mu}(\frkX)u]+[Tv,\Courant{\frkX,Tu}]-T(\rho(T(D_{\rho,\mu}(\frkX)u))v)+T(\rho(\Courant{\frkX,Tu})v)\\
~ &&-T(D_{\rho,\mu}(\frkX)[u,v]_T)+\Courant{\frkX,T[u,v]_T}\\
~ &\stackrel{\eqref{Ooperator1}}{=}&T(\rho(Tu)D_{\rho,\mu}(\frkX)v)-T(\rho(T(D_{\rho,\mu}(\frkX)v))u)-[Tu,\Courant{\frkX,Tv}]+T(\rho(T(D_{\rho,\mu}(\frkX)v))u)-T(\rho(\Courant{\frkX,Tv})u)\\
~ &&-T(\rho(Tv)D_{\rho,\mu}(\frkX)u)+T(\rho(T(D_{\rho,\mu}(\frkX)u))v)+[Tv,\Courant{\frkX,Tu}]-T(\rho(T(D_{\rho,\mu}(\frkX)u))v)+T(\rho(\Courant{\frkX,Tu})v)\\
~ &&-T(D_{\rho,\mu}(\frkX)(\rho(Tu)v-\rho(Tv)u))+\Courant{\frkX,[Tu,Tv]}\\
~ &\stackrel{\eqref{LY3},\eqref{RLYe}}{=}&0,
\end{eqnarray*}}
and}
{\footnotesize
\begin{eqnarray*}
~ &&\delta_{\rm II}^T\Big(\delta(\frkX)\Big)(u,v,w)\\
 ~ &\stackrel{\eqref{2cochain}}{=}&-\delta(\frkX)(\Courant{u,v,w}_T)+D_{\varrho,\varpi}(u,v)\Big(\delta (\frkX)w\Big)+\varpi(v,w)\Big(\delta (\frkX)u\Big)-\varpi(u,w)\Big(\delta (\frkX)v\Big)\\
 ~ &\stackrel{\eqref{repre1},\eqref{repre2}}{=}&\Courant{Tu,Tv,TD_{\rho,\mu}(\frkX)w-\Courant{\frkX,Tw}}+\Courant{TD_{\rho,\mu}(\frkX)u-\Courant{\frkX,Tu},Tv,Tw}\\
 ~ &&+\Courant{Tu,TD_{\rho,\mu}(\frkX)v-\Courant{\frkX,Tv},Tw}\\
 ~ &&-TD(\frkX)\Big(D_{\rho,\mu}(Tu,Tv)w+\mu(Tv,Tw)u-\mu(Tu,Tw)v\Big)\\
 ~ &&+\Courant{\frkX,T\big(D_{\rho,\mu}(Tu,Tv)w+\mu(Tv,Tw)u-\mu(Tu,Tw)v\big)}\\
 ~ &&-T\Big(D_{\rho,\mu}\big(TD_{\rho,\mu}(\frkX)u-\Courant{\frkX,Tu},Tv\big)w-D_{\rho,\mu}\big(TD_{\rho,\mu}(\frkX)v-\Courant{\frkX,Tv},Tu\big)w\Big)\\
 ~ &&-T\Big(\mu\big(TD_{\rho,\mu}(\frkX)v-\Courant{\frkX,Tv},Tw\big)u-\mu\big(TD_{\rho,\mu}(\frkX)u-\Courant{\frkX,Tu},Tw\big)v\Big)\\
 ~ &&-T\Big(\mu\big(Tv,TD_{\rho,\mu}(\frkX)w-\Courant{\frkX,Tw}\big)u-\mu\big(Tu,TD_{\rho,\mu}(\frkX)w-\Courant{\frkX,Tw}\big)v\Big)\\
 ~ &\stackrel{\eqref{Ooperator2}}{=}&\Courant{Tu,Tv,TD_{\rho,\mu}(\frkX)w}-\Courant{Tu,Tv,\Courant{\frkX,Tw}}+\Courant{TD_{\rho,\mu}(\frkX)u,Tv,Tw}\\
 ~ &&-\Courant{\Courant{\frkX,Tu},Tv,Tw}+\Courant{Tu,TD_{\rho,\mu}(\frkX)v,Tw}-\Courant{Tu,\Courant{\frkX,Tv},Tw}\\
 ~ &&-TD_{\rho,\mu}(\frkX)\Big(D_{\rho,\mu}(Tu,Tv)w+\mu(Tv,Tw)u-\mu(Tu,Tw)v\Big)+\Courant{\frkX,\Courant{Tu,Tv,Tw}}\\
 ~ &&-T\Big(D_{\rho,\mu}\big(TD_{\rho,\mu}(\frkX)u-\Courant{\frkX,Tu},Tv\big)w-D_{\rho,\mu}\big(TD_{\rho,\mu}(\frkX)v-\Courant{\frkX,Tv},Tu\big)w\Big)\\
 ~ &&-T\Big(\mu\big(TD_{\rho,\mu}(\frkX)v-\Courant{\frkX,Tv},Tw\big)u-\mu\big(TD_{\rho,\mu}(\frkX)u-\Courant{\frkX,Tu},Tw\big)v\Big)\\
 ~ &&-T\Big(\mu\big(Tv,\big(TD_{\rho,\mu}(\frkX)w-\Courant{\frkX,Tw}\big)u-\mu\big(Tu,TD_{\rho,\mu}(\frkX)w-\Courant{\frkX,Tw}\big)v\Big)\\
 ~ &\stackrel{\eqref{fundamental}}{=}&\Courant{Tu,Tv,TD_{\rho,\mu}(\frkX)w}+\Courant{TD_{\rho,\mu}(\frkX)u,Tv,Tw}+\Courant{Tu,TD_{\rho,\mu}(\frkX)v,Tw}\\
 ~ &&-TD_{\rho,\mu}(\frkX)\Big(D_{\rho,\mu}(Tu,Tv)w+\mu(Tv,Tw)u-\mu(Tu,Tw)v\Big)\\
 ~ &&-T\Big(D_{\rho,\mu}\big(TD_{\rho,\mu}(\frkX)u,Tv\big)w-D_{\rho,\mu}\big(\Courant{\frkX,Tu},Tv\big)w-D_{\rho,\mu}\big(TD_{\rho,\mu}(\frkX)v,Tu\big)w\Big)
 +D_{\rho,\mu}\big(\Courant{\frkX,Tv},Tu\big)w\Big)\\
 ~ &&-T\Big(\mu\big(TD_{\rho,\mu}(\frkX)v,Tw\big)u-\mu\big(\Courant{\frkX,Tv},Tw\big)u-\mu\big(TD_{\rho,\mu}(\frkX)u,Tw\big)v\Big)+\mu\big(\Courant{\frkX,Tu},Tw\big)v\Big)\\
 ~ &&-T\Big(\mu\big(Tv,TD_{\rho,\mu}(\frkX)w\big)u-\mu\big(Tv,\Courant{\frkX,Tw}\big)u-\mu\big(Tu,TD_{\rho,\mu}(\frkX)w\big)v+\mu\big(Tu,\Courant{\frkX,Tw}\big)v\Big)\\
 ~ &\stackrel{\eqref{Ooperator2}}{=}&T\Big(D_{\rho,\mu}(Tu,Tv)D_{\rho,\mu}(\frkX)w+\mu(Tv,TD_{\rho,\mu}(\frkX)w)u-\mu(Tu,TD_{\rho,\mu}(\frkX)w)v\Big)\\
 ~ &&+T\Big(D_{\rho,\mu}(TD_{\rho,\mu}(\frkX)u,Tv)w+\mu(Tv,Tw)D_{\rho,\mu}(\frkX)u-\mu(TD_{\rho,\mu}(\frkX)u,Tw)v\Big)\\
 ~ &&+T\Big(D_{\rho,\mu}(Tu,TD_{\rho,\mu}(\frkX)v)w+\mu(TD_{\rho,\mu}(\frkX)v,Tw)u-\mu(Tu,Tw)D_{\rho,\mu}(\frkX)v\Big)\\
 ~ &&-TD_{\rho,\mu}(\frkX)\Big(D_{\rho,\mu}(Tu,Tv)w+\mu(Tv,Tw)u-\mu(Tu,Tw)v\Big)\\
 ~ &&-T\Big(D_{\rho,\mu}\big(TD_{\rho,\mu}(\frkX)u,Tv\big)w-D_{\rho,\mu}\big(\Courant{\frkX,Tu},Tv\big)w-D_{\rho,\mu}\big(TD_{\rho,\mu}(\frkX)v,Tu\big)w
 +D_{\rho,\mu}\big(\Courant{\frkX,Tv},Tu\big)w\Big)\\
 ~ &&-T\Big(\mu\big(TD_{\rho,\mu}(\frkX)v,Tw\big)u-\mu\big(\Courant{\frkX,Tv},Tw\big)u-\mu\big(TD_{\rho,\mu}(\frkX)u,Tw\big)v)+\mu\big(\Courant{\frkX,Tu},Tw\big)v\Big)\\
 ~ &&-T\Big(\mu\big(Tv,TD_{\rho,\mu}(\frkX)w\big)u-\mu\big(Tv,\Courant{\frkX,Tw}\big)u-\mu\big(Tu,TD_{\rho,\mu}(\frkX)w\big)v+\mu\big(Tu,\Courant{\frkX,Tw}\big)v\Big)\\
 ~ &=&T\Big(D_{\rho,\mu}(Tu,Tv)D_{\rho,\mu}(\frkX)w+\mu(Tv,Tw)D_{\rho,\mu}(\frkX)u-\mu(Tu,Tw)D_{\rho,\mu}(\frkX)v\Big)\\
 ~ &&-TD_{\rho,\mu}(\frkX)\Big(D_{\rho,\mu}(Tu,Tv)w+\mu(Tv,Tw)u-\mu(Tu,Tw)v\Big)\\
 ~ &&+T\Big(D_{\rho,\mu}(\Courant{\frkX,Tu},Tv)w-D_{\rho,\mu}(\Courant{\frkX,Tv},Tu)w+\mu(\Courant{\frkX,Tv},Tw)u\\
 ~ &&-\mu(\Courant{\frkX,Tu},Tw)v+\mu(Tv,\Courant{\frkX,Tw})u-\mu(Tu,\Courant{\frkX,Tw})v\Big)\\
 &\stackrel{\eqref{RLY5},\eqref{RLY5a}}{=}&0.
\end{eqnarray*}}
Similarly, we deduce that $\delta^T_{\rm I}\Big(\delta(\frkX)\Big)(u,v)=0$ for all $u,v\in V$.
This finishes the proof.
\end{proof}

Let $T$ be a relative Rota-Baxter operator on a Lie-Yamaguti algebra $(\g,[\cdot,\cdot],\Courant{\cdot,\cdot,\cdot})$ with  respect to a representation $(V;\rho,\mu)$. Combining the Yamaguti cohomology with Proposition \ref{0cocy}, we obtain a well-defined cochain complex $(\huaC_T^\bullet(V,\g)=\bigoplus_{n=0}^\infty\huaC_T^n(V,\g),\de)$, where the $n$-cochians $\huaC^n_T(V,\g)$ and the coboundary map $\de:\huaC^n_{T}(V,\g)\to \huaC^{n+1}_{T}(V,\g)$ are defined to be
\begin{eqnarray*}
\huaC^n_T(V,\g):=
\begin{cases}
C^n_{\rm LieY}(V,\g),&n\geqslant 1,\\
\wedge^2\g,&n=0,
\end{cases}
\end{eqnarray*}
and
\begin{eqnarray*}
\de:=
\begin{cases}
\delta^T,&n\geqslant 1,\\
\delta,&n=0,
\end{cases}
\end{eqnarray*}
respectively. Here $\delta^T:C_{\rm LieY}^{n}(V,\g)\to C_{\rm LieY}^{n+1}(V,\g)~(n\geqslant 1)$ is the corresponding Yamaguti coboundary operator on the sub-adjacent Lie-Yamaguti algebra $(V,[\cdot,\cdot]_T,\Courant{\cdot,\cdot,\cdot}_T)$ with coefficients in the representation $(\g;\varrho,\varpi)$.
\begin{defi}\label{cohomology}
The cohomology of cochian complex $(\huaC_T^\bullet(V,\g)=\bigoplus_{n=0}^\infty\huaC_T^n(V,\g),\de)$ is called the {\bf cohomology of relative Rota-Baxter operator $T$} on Lie-Yamaguti algebra $(\g,[\cdot,\cdot],\Courant{\cdot,\cdot,\cdot})$ with respect to the representation $(V;\rho,\mu)$. Denote the set of $n$-cocycles and $n$-coboundaries by $\huaZ^n_T(V,\g)$ and $\huaB^n_T(V,\g)$ respectively. The $n$-th cohomology group of relative Rota-Baxter operator $T$ is taken to be
\begin{eqnarray*}
\huaH^n_T(V,\g):=\huaZ^n_T(V,\g)/\huaB^n_T(V,\g), \quad n\geqslant 1.
\end{eqnarray*}
\end{defi}

Let us give the formula of $\delta^T$ explicitly:

\begin{itemize}
\item if $n\geqslant 1$, $\delta^T:C_{\rm LieY}^{n+1}(V,\g)\to C_{\rm LieY}^{n+2}(V,\g)$ is given by for any $F=(f,g)\in C_{\rm LieY}^{n+1}(V,\g)$,
{\footnotesize
\begin{eqnarray}
~\nonumber &&\Big(\delta^T_{\rm I}(F)\Big)(\huaV_1,\cdots,\huaV_{n+1})\\
~\label{Ocohomo1} &=&(-1)^{n}\Big([Tu_{n+1},g(\huaV_1,\cdots,\huaV_n,v_{n+1})]-[Tv_{n+1},g(\huaV_1,\cdots,\huaV_n,u_{n+1})]\\
~\nonumber &&-g(\huaV_1,\cdots,\huaV_n,\rho(Tu_{n+1})v_{n+1}-\rho(Tv_{n+1})u_{n+1})+T\big(\rho(g(\huaV_1,\cdots,\huaV_n,v_{n+1}))u_{n
+1}\big)\\
~\nonumber &&-T\big(\rho(g(\huaV_1,\cdots,\huaV_n,u_{n+1}))v_{n
+1}\big)\Big)\\
~ &&\nonumber+\sum_{k=1}^{n+1}(-1)^{k+1}\Big(\Courant{Tu_{k},Tv_k,f(\huaV_1,\cdots,\hat{\huaV_k},\cdots,\huaV_{n+1})}+T\big(\mu(Tv_k,f(\huaV_1,\cdots,\hat{\huaV_k},\cdots,\huaV_{n+1}))u_k\big)\\
~ \nonumber&&-T\big(\mu(Tu_k,f(\huaV_1,\cdots,\hat{\huaV_k},\cdots,\huaV_{n+1}))v_k\big)\Big)\\
~ \nonumber&&+\sum_{1\leqslant k<l\leqslant n+1}(-1)^kf(\huaV_1,\cdots,\hat{\huaV_k},\cdots,\huaV_k\circ\huaV_l,\cdots,\huaV_{n+1}),
\end{eqnarray}}
and
{\footnotesize
\begin{eqnarray}
~\nonumber &&\Big(\delta^T_{\rm II}(F)\Big)(\huaV_1,\cdots,\huaV_{n+1},w)\\
~\label{Ocohomo2}&=&(-1)^n\Big(\Courant{g(\huaV_1,\cdots,\huaV_n,u_{n+1}),Tv_{n+1},Tw}-\Courant{g(\huaV_1,\cdots,\huaV_n,v_{n+1}),Tu_{n+1},Tw}\\
~\nonumber &&+T\big(D(g(\huaV_1,\cdots,\huaV_n,u_{n+1}),Tv_{n+1})w-\mu(g(\huaV_1,\cdots,\huaV_n,u_{n+1}),Tw)v_{n+1}\\
~\nonumber &&+D(g(\huaV_1,\cdots,\huaV_n,v_{n+1}),Tu_{n+1})w-\mu(g(\huaV_1,\cdots,\huaV_n,v_{n+1}),Tw)u_{n+1}\big)\Big)\\
~ \nonumber&&+\sum_{k=1}^{n+1}(-1)^{k}\Big(\Courant{Tu_k,Tv_k,g(\huaV_1,\cdots,\hat{\huaV_k},\cdots,\huaV_{n+1},w)}+T\big(\mu(Tv_k,g(\huaV_1,\cdots,\hat{\huaV_k},\cdots,\huaV_{n+1},w))u_k\\
~\nonumber &&-\mu(Tu_k,g(\huaV_1,\cdots,\hat{\huaV_k},\cdots,\huaV_{n+1},w))v_k\big)\Big)\\
~\nonumber &&+\sum_{1\leqslant k<l\leqslant n+1}(-1)^kg(\huaV_1,\cdots,\hat{\huaV_k},\cdots,\huaV_k\circ\huaV_l,\cdots,\huaV_{n+1},w)\\
~\nonumber &&+\sum_{k=1}^{n+1}(-1)^kg(\huaV_1,\cdots,\hat{\huaV_k},\cdots,\huaV_{n+1},\Courant{u_k,v_k,w}_T),
\end{eqnarray}}
where $\huaV_i=u_i\wedge v_i\in \wedge^2V,~ (1\leqslant i\leqslant n+1),~ w\in V$ and $\huaV_k\circ\huaV_l=\Courant{u_k,v_k,u_l}_T\wedge v_l+u_l\wedge\Courant{u_k,v_k,v_l}_T$.

\item if $n=0$, for any $f\in C_{\rm LieY}^1(V,\g)=\Hom(V,\g)$,
$$\delta^T:C_{\rm LieY}^1(V,\g)\to C_{\rm LieY}^2(V,\g),\quad f \mapsto (\delta^T_{\rm I}(f),\delta^T_{\rm II}(f))$$
is given by
\begin{eqnarray*}
(\delta^T_{\rm I}(f))(u,v)&=&[Tu,f(v)]-[Tv,f(u)]+T\Big(\rho(f(v)u)-\rho(f(u)v)\Big)-f([u,v]_T),\\
(\delta^T_{\rm II}(f))(u,v,w)&=&\Courant{Tu,Tv,f(w)}+\Courant{f(u),Tv,Tw}-\Courant{f(v),Tu,Tw}-f(\Courant{u,v,w}_T)\\
~ &&-T\Big(D(f(u),Tv)w-D(f(v),Tu)w+\mu(Tv,f(w))u-\mu(Tu,f(w))v\\
~ &&-\mu(f(u),Tw)v+\mu(f(v),Tw)u\Big), \qquad \forall u,v,w \in V.
\end{eqnarray*}
\end{itemize}

At the end of this section, we show that the coboundary map $\de$ coincides (up to sign) with the differential $l_1^T$ defined by
Equation \eqref{diffs} which involves the Maurer-Cartan element $T$ of the twisted $L_\infty$-algebra $(C^*(V,\g),l_1^T,l_2^T,l_3^T)$.

\begin{thm}\label{diff}
Let $T:V\to \g$ be a relative Rota-Baxter operator on a Lie-Yamaguti algebra $(\g,[\cdot,\cdot],\Courant{\cdot,\cdot,\cdot})$ with respect to a representation $(V;\rho,\mu)$. Then we have
\begin{eqnarray}
\label{Oo1}\de(F)&=&(-1)^{n-1}l_1^T(F),\quad \forall F\in C^n(V,\g),~n=1,2,3,\cdots,\\
\label{Oo2}\de(f)&=&l_1^T(f),\quad \forall f\in C^0(V,\g)=\Hom(V,\g).
\end{eqnarray}
\end{thm}
\begin{proof}
First, for all $x,y, z\in \g$, $u,v,w \in V$, we compute that
\begin{eqnarray*}
~ &&[\Pi+\Theta,T]_{\rm I}(x+u,y+v)\\
~ &=&\Big((\Pi+\Theta)\circ T\Big)_{\rm I}(x+u,y+v)-\Big( T\circ (\Pi+\Theta)\Big)_{\rm I}(x+u,y+v)\\
~ &=&[Tu,y]+\rho(Tu)v+[x,Tv]-\rho(Tv)u-T\Big(\rho(x)v-\rho(y)u\Big),\\
~ &&\\
~ &&[\Pi+\Theta,T]_{\rm II}(x+u,y+v,z+w)\\
~ &=&\Big((\Pi+\Theta)\circ T\Big)_{\rm II}(x+u,y+v,z+w)-\Big( T\circ (\Pi+\Theta)\Big)_{\rm II}(x+u,y+v,z+w)\\
~ &=&\Courant{Tu,y,z}+D(Tu,y)w-\mu(Tu,z)v+\Courant{x,Tv,z}+D(x,Tv)w+\mu(Tv,z)u\\
~ &&+\Courant{x,y,Tw}+\mu(y,Tw)u-\mu(x,Tw)v-T\Big(D(x,y)w+\mu(y,z)u-\mu(x,z)v\Big),\\
~ &&\\
~ &&[[\Pi+\Theta,T]_{\mathsf{LieY}},T]_{\rm I}(x+u,y+v)\\
~ &=&\Big([\Pi+\Theta,T]_{\mathsf{LieY}}\circ T\Big)_{\rm I}(x+u,y+v)-\Big( T\circ [\Pi+\Theta,T]_{\mathsf{LieY}}\Big)_{\rm I}(x+u,y+v)\\
~ &=&2\Big([Tu,Tv]-T\big(\rho(Tu)v-\rho(Tv)u\big)\Big),\\
~ &&\\
~ &&[[\Pi+\Theta,T]_{\mathsf{LieY}},T]_{\rm II}(x+u,y+v,z+w)\\
~ &=&[\Pi+\Theta,T]_{\rm II}(Tu,y+v,z+w)+[\Pi+\Theta,T]_{\rm II}(x+u,Tv,z+w)\\
~ &&+[\Pi+\Theta,T]_{\rm II}(x+u,y+v,Tw)-T\Big([\Pi+\Theta,T]_{\rm II}(x+u,y+v,z+w)\Big)\\
~ &=&2\Big(\Courant{Tu,Tv,z}+\Courant{x,Tv,Tw}+\Courant{Tu,y,Tw}+D(Tu,Tv)w+\mu(Tv,Tw)u-\mu(Tu,Tw)v\Big)\\
~ &&-2T\Big(D(Tu,y)w+D(x,Tv)w-\mu(Tu,z)v-\mu(x,Tw)v+\mu(Tv,z)u+\mu(y,Tw)u\Big).
\end{eqnarray*}
Then, by \eqref{multi1} and \eqref{multi2}, for all $F=(f,g)\in C^n(V,\g)~(n\geqslant 1)$, let us compute that
{\footnotesize
\begin{eqnarray*}
~ &&(l_2)_{\rm I}\Big(T,F\Big)(\huaV_1,\cdots,\huaV_{n+1})\\
~ &=&\pr_\g[[\Pi+\Theta,T]_{\mathsf{LieY}},F]_{\rm I}(\huaV_1,\cdots,\huaV_{n+1})\\
~ &=&\Big([\Pi+\Theta,T]_{\mathsf{LieY}}\circ F\Big)_{\rm I}(\huaV_1,\cdots,\huaV_{n+1})-(-1)^n\Big(F\circ [\Pi+\Theta,T]_{\mathsf{LieY}}\Big)_{\rm I}(\huaV_1,\cdots,\huaV_{n+1})\\
~ &=&\sum_{\sigma\in\mathbb S_{(1,n)}\atop\sigma(n+1)=n+1}(-1)^n sign(\sigma)[\Pi+\Theta,T]_{\rm II}(\huaV_{\sigma(1)},f(\huaV_{\sigma(2)},\cdots,\huaV_{\sigma(n+1)}))\\
~ &&+[\Pi+\Theta,T]_{\rm I}(u_{n+1},g(\huaV_1,\cdots,\huaV_n,v_{n+1}))+[\Pi+\Theta,T]_{\rm I}(g(\huaV_1,\cdots,\huaV_n,u_{n+1}),v_{n+1})\\
~ &&-(-1)^n\Big((-1)^ng(\huaV_1,\cdots,\huaV_n,[\Pi+\Theta,T]_{\rm I}(\huaV_{n+1}))\\
~ &&+\sum_{k=1}^n(-1)^{k-1}\sum_{\sigma\in\mathbb S_{(k-1,1)}}sign(\sigma)f(\huaV_{\sigma(1)},\cdots,\huaV_{\sigma(k-1)},u_{k+1}\wedge [\Pi+\Theta,T]_{\rm II}(\huaV_{\sigma(k)},v_{k+1}),\huaV_{k+2},\cdots,\huaV_{n+1})\\
~ &&+\sum_{k=1}^n(-1)^{k-1}\sum_{\sigma\in\mathbb S_{(k-1,1)}}sign(\sigma)f(\huaV_{\sigma(1)},\cdots,\huaV_{\sigma(k-1)}, [\Pi+\Theta,T]_{\rm II}(\huaV_{\sigma(k)},u_{k+1})\wedge v_{k+1},\huaV_{k+2},\cdots,\huaV_{n+1})\Big)\\
~ &=&[Tu_{n+1},g(\huaV_1,\cdots,\huaV_n,v_{n+1})]+T(\rho(g(\huaV_1,\cdots,\huaV_n,v_{n+1}))u_{n+1})+[g(\huaV_1,\cdots,\huaV_n,u_{n+1}),Tv_{n+1}]\\
~ &&-T(\rho(g(\huaV_1,\cdots,\huaV_n,u_{n+1}))v_{n+1})-g(\huaV_1,\cdots,\huaV_n,\rho(Tu_{n+1})v_{n+1}-\rho(Tv_{n+1})u_{n+1}),
\end{eqnarray*}}
and
{\footnotesize
\begin{eqnarray*}
~ &&(l_3)_{\rm I}\Big(T,T,F\Big)(\huaV_1,\cdots,\huaV_{n+1})\\
~ &=&\pr_\g[[[\Pi+\Theta,T]_{\mathsf{LieY}}],T]_{\mathsf{LieY}},F]_{\rm I}(\huaV_1,\cdots,\huaV_{n+1})\\
~ &=&\Big([[\Pi+\Theta,T]_{\mathsf{LieY}},T]_{\mathsf{LieY}}\circ F\Big)_{\rm I}(\huaV_1,\cdots,\huaV_{n+1})-(-1)^n\Big(F\circ [[\Pi+\Theta,T]_{\mathsf{LieY}},T]_{\mathsf{LieY}}\Big)_{\rm I}(\huaV_1,\cdots,\huaV_{n+1})\\
~ &=&\sum_{\sigma\in\mathbb S_{(1,n)}\atop \sigma(n+1)=n+1}(-1)^nsign(\sigma)[[\Pi+\Theta,T]_{\mathsf{LieY}},T]_{\rm II}(\huaV_{\sigma(1)},f(\huaV_{\sigma(2)},\cdots,\huaV_{\sigma(n+1)}))\\
~&&+[[\Pi+\Theta,T]_{\mathsf{LieY}},T]_{\rm I}(u_{n+1},g(\huaV_{1},\cdots,\huaV_n,v_{n+1}))+[[\Pi+\Theta,T]_{\mathsf{LieY}},T]_{\rm I}(g(\huaV_{1},\cdots,\huaV_n,u_{n+1}),v_{n+1})\\
~ &&-(-1)^n\Big(\sum_{k=1}^n(-1)^{k-1}\sum_{\sigma\in\mathbb S_{(k-1,1)}}sign(\sigma)f(\huaV_{\sigma(1)},\cdots,\huaV_{\sigma(k-1)},u_{k+1}\wedge[[\Pi+\Theta,T]_{\mathsf{LieY}},T]_{\rm II}(\huaV_{\sigma(k)},v_{k+1}),\huaV_{k+2},\cdots,\huaV_{n+1})\\
~ &&+\sum_{k=1}^n(-1)^{k-1}\sum_{\sigma\in\mathbb S_{(k-1,1)}}sign(\sigma)f(\huaV_{\sigma(1)},\cdots,\huaV_{\sigma(k-1)},[[\Pi+\Theta,T]_{\mathsf{LieY}},T]_{\rm II}(\huaV_{\sigma(k)},u_{k+1})\wedge v_{k+1},\huaV_{k+2},\cdots,\huaV_{n+1})\Big)\\
~ &=&2\Big(\sum_{\sigma\in\mathbb S_{(1,n)}\atop \sigma(n+1)=n+1}(-1)^nsign(\sigma)\big(\Courant{Tu_{\sigma(1)},Tv_{\sigma(1)},f(\huaV_{\sigma(2)},\cdots,\huaV_{\sigma(n+1)})}+T(\mu(Tv_{\sigma(1)},
f(\huaV_{\sigma(2)},\cdots,\huaV_{\sigma(n+1)}))u_{\sigma(1)})\\
~ &&-T(\mu(Tu_{\sigma(1)},
f(\huaV_{\sigma(2)},\cdots,\huaV_{\sigma(n+1)}))v_{\sigma(1)})\big)+\sum_{k=1}^n(-1)^{k-1}\sum_{\sigma\in \mathbb S_{(k-1,1)}}f(\huaV_{\sigma(1)},\cdots,\huaV_{\sigma(k-1)},\huaV_{\sigma(k)}\circ \huaV_{k+1},\huaV_{k+2},\cdots,\huaV_{n+1})\Big).
\end{eqnarray*}}
Thus, we have that
\begin{eqnarray*}
~ &&(-1)^{n-1}\Big(l_1^T\Big)_{\rm I}(F)(\huaV_1,\cdots,\huaV_{n+1})\\
~ &\stackrel{\eqref{diffs}}{=}&(-1)^{n-1}\Big((l_2)_{\rm I}(T,F)+\half (l_3)_{\rm I}(T,T,F)\Big)(\huaV_1,\cdots,\huaV_{n+1})\\
~ &\stackrel{\eqref{multi1},\eqref{multi2}}{=}&(-1)^{n}\Big([Tu_{n+1},g(\huaV_1,\cdots,\huaV_n,v_{n+1})]-[Tv_{n+1},g(\huaV_1,\cdots,\huaV_n,u_{n+1})]\\
~ &&-g(\huaV_1,\cdots,\huaV_n,\rho(Tu_{n+1})v_{n+1}-\rho(Tv_{n+1})u_{n+1})+T\big(\rho(g(\huaV_1,\cdots,\huaV_n,v_{n+1}))u_{n
+1}\big)\\
~ &&-T\big(\rho(g(\huaV_1,\cdots,\huaV_n,u_{n+1}))v_{n
+1}\big)\Big)\\
~ &&+\sum_{k=1}^{n+1}(-1)^{k+1}\Big(\Courant{Tu_{k},Tv_k,f(\huaV_1,\cdots,\hat{\huaV_k},\cdots,\huaV_{n+1})}+T\big(\mu(Tv_k,f(\huaV_1,\cdots,\hat{\huaV_k},\cdots,\huaV_{n+1}))u_k\big)\\
~ &&-T\big(\mu(Tu_k,f(\huaV_1,\cdots,\hat{\huaV_k},\cdots,\huaV_{n+1}))v_k\big)\Big)\\
~ &&+\sum_{1\leqslant k<l\leqslant n+1}(-1)^kf(\huaV_1,\cdots,\hat{\huaV_k},\cdots,\huaV_k\circ\huaV_l,\cdots,\huaV_{n+1}),\\
~ &\stackrel{\eqref{Ocohomo1}}{=}&\delta_{\rm I}^T(F)(\huaV_1,\cdots,\huaV_{n+1}).
\end{eqnarray*}
Hence, we have
\begin{eqnarray}
\delta_{\rm I}^T(F)=(-1)^{n-1}\Big(l_1^T\Big)_{\rm I}(F), \quad\forall F\in C^n(V,\g) ~(n\geqslant 1).\label{Oopcoho1}
\end{eqnarray}
Similarly, by a direct computation, we have that for all $F\in C^n(V,\g) ~(n\geqslant 1)$ and for all $\huaV_i\in \wedge^2V~(i=1,\cdots,n+1)$ and $u\in V$
 \begin{eqnarray*}
 ~(l_2)_{\rm II}(F)(\huaV_1,\cdots,\huaV_{n+1},u)&=&\pr_\g[[\Pi+\Theta,T]_{\mathsf{LieY}},F]_{\rm II}(\huaV_1,\cdots,\huaV_{n+1},u)=0,\\
 (l_3)_{\rm II}(T,T,F)(\huaV_1,\cdots,\huaV_{n+1},u)&=&(-1)^{n-1}2\delta_{\rm II}^T(F)(\huaV_1,\cdots,\huaV_{n+1},u).
 \end{eqnarray*}
Thus, we have that
\begin{eqnarray}
(-1)^{n-1}\Big(l_1^T\Big)_{\rm II}(F)(\huaV_1,\cdots,\huaV_{n+1},u)=\delta_{\rm II}^T(F)(\huaV_1,\cdots,\huaV_{n+1},u).\label{Oopcoho2}
\end{eqnarray}
Combining Eqs. \eqref{Oopcoho1} and \eqref{Oopcoho2} yields \eqref{Oo1}. Moreover, for any $f\in \Hom(V,\g)$, by a direct computation, we have that for all $u,v,w\in V$,
\begin{eqnarray*}
~ &&\Big(l_2\Big)_{\rm I}(T,f)(u,v)\\
~ &=&\pr_\g[[\Pi+\Theta,T]_{\mathsf{LieY}},f]_{\rm I}(u,v)\\
~ &=&[f(u),Tv]+[Tu,f(v)]-T\Big(\rho(f(u))v-\rho(f(v))u\Big)-f([u,v]_T),\\
~ &&\Big(l_3\Big)_{\rm I}(T,T,f)(u,v)=\pr_\g[[[\Pi+\Theta,T]_{\mathsf{LieY}},T]_{\mathsf{LieY}},f]_{\rm I}=0.
\end{eqnarray*}
Thus we get
\begin{eqnarray}
\Big(l_1^T\Big)_{\rm I}(f)(u,v)=\Big(l_2\Big)_{\rm I}(T,f)(u,v)=\delta_{\rm I}^T(f)(u,v). \label{Oopcoho3}
\end{eqnarray}
Similarly,
\begin{eqnarray*}
~ &&\Big(l_2\Big)_{\rm II}(T,f)(u,v,w)=\pr_\g[[\Pi+\Theta,T]_{\mathsf{LieY}},f]_{\rm II}(u,v,w)=0,\\
~ &&\Big(l_3\Big)_{\rm II}(T,T,f)(u,v,w)=\pr_\g[[[\Pi+\Theta,T]_{\mathsf{LieY}},T]_{\mathsf{LieY}},f]_{\rm II}(u,v,w)\\
~ &=&2\Big(\Courant{f(u),Tv,Tw}-T\big(D(f(u),Tv)w-\mu(f(u),Tw)v+c.p.\big)-f(\Courant{u,v,w}_T)\Big),
\end{eqnarray*}
where the last equation is equivalent to $\Big(l_3\Big)_{\rm II}(T,T,f)(u,v,w)=2\delta_{\rm II}^T(f)(u,v,w)$. Thus, for all $u,v,w\in V$, we have that
\begin{eqnarray}
\Big(l_1^T\Big)_{\rm II}(f)(u,v,w)=\delta_{\rm II}^T(f)(u,v,w). \label{Oopcoho4}
\end{eqnarray}
Hence, Eqs. \eqref{Oopcoho3} and \eqref{Oopcoho4} give \eqref{Oo2}.
This completes the proof.
\end{proof}

\section{Deformations of relative Rota-Baxter operators}

After establishing cohomology theory of relative Rota-Baxter operators, we study deformations of relative Rota-Baxter operators, i.e., we use this cohomology to characterize deformations of relative Rota-Baxter operators on Lie-Yamaguti algebras.

\subsection{Linear deformations of relative Rota-Baxter operators}

In this subsection, we explore linear deformations of relative Rota-Baxter operators on Lie-Yamaguti algebras, and we show that the infinitesimals of two equivalent linear deformations of a relative Rota-Baxter operator on Lie-Yamaguti algebra are in the same cohomology classes of the first cohomology group.
\begin{defi}
Let $T$ and $T'$ be two relative Rota-Baxter operators on a Lie-Yamaguti algebra $(\g,[\cdot,\cdot],\Courant{\cdot,\cdot,\cdot})$ with respect to a representation $(V;\rho,\mu)$. A {\bf homomorphism} from $T'$ to $T$ consists of a Lie-Yamaguti homomorphism $\phi_\g: \g \to \g$ and a linear map $\phi_V: V \to V$ such that
\begin{eqnarray}
\label{homo1}T\circ \phi_V&=&\phi_\g\circ T'\\
~\label{homo2}\phi_V(\rho(x)v)&=&\rho(\phi_\g(x))\phi_V(v),\\
~\label{homo3}\phi_V\mu(x,y)(v)&=&\mu(\phi_\g(x),\phi_\g(y))(\phi_V(v)), \quad \forall x,y \in \g,~v\in V.
\end{eqnarray}
In particular, if $\phi_\g$ and $\phi_V$ are invertible, then $(\phi_\g,\phi_V)$ is called an {\bf isomorphism} from $T'$ to $T$.
\end{defi}

\begin{rmk}
If two relative Rota-Baxter operators $T$ and $T'$ are homomorphic, that is, there exists a pair $(\phi_\g,\phi_V)$ such that \eqref{homo1}-\eqref{homo3} hold, then by a direct computation we have that
\begin{eqnarray*}
\label{homo4}\phi_VD(x,y)(v)&=&D(\phi_\g(x),\phi_\g(y))(\phi_V(v)), \quad \forall x,y \in \g,~v\in V.
\end{eqnarray*}
\end{rmk}
\emptycomment{
\begin{pro}
Let $T$ and $T'$ be two relative Rota-Baxter operators on a Lie-Yaamguti algebra $(\g,[\cdot,\cdot],\Courant{\cdot,\cdot,\cdot})$ with respect to a representation $(V;\rho,\mu)$ and $(\phi_\g,\phi_V)$ a homomorphism from $T'$ to $T$. Then $\phi_V$ is a homomorphism from pre-Lie-Yamaguti algebra\footnote{For the notion of pre-Lie-Yamaguti algebras, one can see \cite{SZ1}.} from $(V,*_{T'},\{\cdot,\cdot,\cdot\}_{T'})$ to $(V,*_T,\{\cdot,\cdot,\cdot\}_{T})$.
\end{pro}
\begin{proof}
For all $u,v,w\in V$, we have
\begin{eqnarray*}
\phi_V(u*_{T'}v)&=&\phi_V(\rho(T'u)v)=\rho(\phi_\g(T'u)\phi_V(v))\\
~ &=&\rho(T(\phi_V(u))\phi_V(v))=\phi_V(u)*_T\phi_V(v),\\
\phi_V(\{u,v,w\}_{T'})&=&\phi_V(\mu(T'v,T'w)u)=\mu(\phi_\g(T'v),\phi_\g(T'w))(\phi_V(u))\\
~ &=&\mu(T(\phi_V(v)),T(\phi_V(w)))(\phi_V(u))=\{\phi_V(u),\phi_V(v),\phi_V(w)\}_{T}.
\end{eqnarray*}
This finishes the proof.
\end{proof}}

\begin{defi}
Let $T$ be a relative Rota-Baxter operator on a Lie-Yamaguti algebra $(\g,[\cdot,\cdot],\Courant{\cdot,\cdot,\cdot})$ with respect to a representation $(V;\rho,\mu)$ and $\frkT:V \to \g$ a linear map. If $T_t=T+t\frkT$ are still relative Rota-Baxter operators on $(\g,[\cdot,\cdot],\Courant{\cdot,\cdot,\cdot})$ for all $t$, we say that $\frkT$ generates a {\bf linear deformation} of the relative Rota-Baxter operator $T$.
\end{defi}
\emptycomment{
It is easy to see that $\frkT$ generates a linear deformation of the relative Rota-Baxter operator $T$ if and only if
\begin{eqnarray}
~\label{cocy1}~ &&[\frkT u,Tv]+[Tu,\frkT v]=T\big(\rho(\frkT u)v-\rho(\frkT v)u\big)+\frkT\big(\rho(Tu)v-\rho(Tv)u\big),\\
~\label{Ooper1}&&[\frkT u,\frkT v]=\frkT \big(\rho(\frkT u)v-\rho(\frkT v)u\big),\\
~ \label{cocy}&&\Courant{Tu,Tv,\frkT w}+\Courant{Tu,\frkT v,Tw}+\Courant{\frkT u,Tv,Tw}\\
~ \nonumber&=&\frkT\Big(D(Tu,Tv)w+\mu(Tv,Tw)v-\mu(Tu,Tw)v\Big)\\
~ &&\nonumber+T\Big(D(Tu,\frkT v)w+D(\frkT u,Tv)w+\mu(Tv,\frkT w)u\\
~ &&\nonumber+\mu(\frkT v,Tw)u-\mu(Tu,\frkT w)v-\mu(\frkT u,Tw)v\Big),\\
~ &&\Courant{\frkT u,\frkT v, Tw}+\Courant{Tu,\frkT v,\frkT w}+\Courant{\frkT u,Tv,\frkT w}\\
~ \nonumber&=&T\Big(D(\frkI u,\frkI v)w+\mu(\frkI v,\frkI w)v-\mu(\frkI u,\frkI w)v\Big)\\
~ &&\nonumber+\frkT\Big(D(Tu,\frkT v)w+D(\frkT u,Tv)w+\mu(Tv,\frkT w)u\\
~ &&\nonumber+\mu(\frkT v,Tw)u-\mu(Tu,\frkT w)v-\mu(\frkT u,Tw)v\Big),\\
~ \label{Oopera}&&\Courant{\frkT u,\frkT v, \frkT w}=\frkT \Big(D(\frkT u,\frkT v)w+\mu(\frkT v,\frkT w)u-\mu(\frkT u,\frkT w)v\Big).
\end{eqnarray}}

\begin{rmk}
If $\frkT$ generates a linear deformation of a relative Rota-Baxter operator $T$, then we have that
\begin{itemize}
\item[(i)] $\frkT$ is a relative Rota-Baxter operator on the Lie-Yamaguti algebra $(\g,[\cdot,\cdot],\Courant{\cdot,\cdot,\cdot})$ with respect to the representation $(V;\rho,\mu)$;
    \item[(ii)]  $\frkT\in\huaC^1(V,\g)$ is a $1$-cocycle of $\delta^T$.
    \end{itemize}
\end{rmk}

\begin{defi}
Let $T:V\to \g$ be a relative Rota-Baxter operator on a Lie-Yamaguti algebra $(\g,[\cdot,\cdot],\Courant{\cdot,\cdot,\cdot})$ with respect to a representation $(V;\rho,\mu)$.
\begin{itemize}
\item [\rm (i)]  Two linear deformations $T_t^1=T+t\frkT_1$ and $T_t^2=T+t\frkT_2$ are said to be {\bf equivalent} if there exists an element $\frkX \in \wedge^2\g$ such that $({\Id}_\g+t\frkL_{\frkX},{\Id}_V+tD(\frkX))$ is a homomorphism from $T_t^2$ to $T_t^1$.
\item [\rm (ii)] A linear deformation $T_t=T+t\frkT$ of a relative Rota-Baxter operator $T$ is said to be {\bf trivial} if there exists an element $\frkX \in \wedge^2\g$ such that $({\Id}_\g+t\frkL_{\frkX},{\Id}_V+tD(\frkX))$ is a homomorphism from $T_t$ to $T$.
    \end{itemize}
\end{defi}

Let $({\Id}_\g+t\frkL_{\frkX},{\Id}_V+tD(\frkX))$ be a homomorphism from $T_t^2$ to $T_t^1$. Then ${\Id}_\g+t\frkL_{\frkX}$ is a Lie-Yamaguti algebra homomorphism of $\g$, which is equivalent to the following conditions:
\begin{eqnarray}
~\label{Nije} [\Courant{\frkX,x},\Courant{\frkX,y}]&=&0,\\
~ \label{Nij1}\Courant{\Courant{\frkX,x},\Courant{\frkX,y},z}+\Courant{\Courant{\frkX,x},y,\Courant{\frkX,z}}+\Courant{x,\Courant{\frkX,y},\Courant{\frkX,z}}&=&0,\\
~ \label{Nij2}\Courant{\Courant{\frkX,x},\Courant{\frkX,y},\Courant{\frkX,z}}&=&0.
\end{eqnarray}

By $T_t^1\big(({\Id_V}+tD(\frkX))v\big)=\big({\Id}_\g+t\frkL_{\frkX}\big)T_t^2(v)$, we have
\begin{eqnarray}
\label{cocycle}(\frkT_2-\frkT_1)(v)&=&T\Big(D(\frkX)v\Big)-\Courant{\frkX,Tv},\\
\frkT_1\Big(D(\frkX)v\Big)&=&\Courant{\frkX,\frkT_2(v)}
\end{eqnarray}
\emptycomment{
By $\Big({\Id}_V+tD(\frkX)\Big)(\rho(x)v)=\rho\Big(({\Id}_\g+t\frkL_\frkX)(x)\Big)({\Id}_V+tD(\frkX))(v)$, we have
\begin{eqnarray}
\rho(\Courant{\frkX,x})D(\frkX)=0.
\end{eqnarray}}

Finally, by $\Big({\Id}_V+tD(\frkX)\Big)\mu(z,w)v=\mu\Big(({\Id}_\g+t\frkL_{\frkX})z,({\Id}_\g+t\frkL_{\frkX})w\Big)({\Id}_V+tD(\frkX))v$, we have
\begin{eqnarray}
\label{Nij3}\mu(z,\Courant{\frkX,w})D(\frkX)+\mu(\Courant{\frkX,z},w)D(\frkX)+\mu(\Courant{\frkX,z},\Courant{\frkX,w})&=&0,\\
~\label{Nij4}\mu(\Courant{\frkX,z},\Courant{\frkX,w})D(\frkX)&=&0.
\end{eqnarray}

Note that \eqref{cocycle} means that there exists $\frkX\in \wedge^2\g$, such that $\frkT_2-\frkT_1=\delta(\frkX)$. Thus we have the following
\begin{thm}\label{thm1}
Let $T:V\to \g$ be a relative Rota-Baxter operator on a Lie-Yamaguti algebra $(\g,[\cdot,\cdot],\Courant{\cdot,\cdot,\cdot})$ with respect to a representation $(V;\rho,\mu)$. If two linear deformations $T_t^1=T+t\frkT_1$ and $T_t^2=T+t\frkT_2$ of $T$ are equivalent, then $\frkT_1$ and $\frkT_2$ are in the same class of the cohomology group $\huaH^1_T(V,\g)$.
\end{thm}

\begin{defi}
Let $T:V\to \g$ be a relative Rota-Baxter operator on a Lie-Yamaguti algebra $(\g,[\cdot,\cdot],\Courant{\cdot,\cdot,\cdot})$ with respect to a representation $(V;\rho,\mu)$. An element $\frkX\in \wedge^2\g$ is called a {\bf Nijenhuis element} with respect to $T$ if $\frkX$ satisfies \eqref{Nije}-\eqref{Nij2}, \eqref{Nij3}, \eqref{Nij4} and the following equation
\begin{eqnarray*}
\Courant{\frkX,T(D(\frkX)v)-\Courant{\frkX,Tv}}=0, \quad \forall v \in V. \label{Nij5}
\end{eqnarray*}
\end{defi}

It is obvious that a trivial deformation of a relative Rota-Baxter operator on a Lie-Yamaguti algebra gives rise to a Nijenhuis element. Indeed, the converse is also true.
\begin{pro}\label{Nijenhuis}
Let $T:V\to \g$ be a relative Rota-Baxter operator on a Lie-Yamaguti algebra $(\g,[\cdot,\cdot],\Courant{\cdot,\cdot,\cdot})$ with respect to a representation $(V;\rho,\mu)$. Then for any Nijenhuis element $\frkX\in \wedge^2 \g$, $T_t:=T+t\frkT$ with $\frkT:=\delta(\frkX)$ is a trivial linear deformation of the relative Rota-Baxter operator $T$.
\end{pro}
\begin{proof}
The proof is similar to the case of Lie algebras or Leibniz algebras etc. Thus we omit the details.
\end{proof}

 \begin{ex}
Let $(\g,[\cdot,\cdot],\Courant{\cdot,\cdot,\cdot})$ be a 2-dimensional Lie-Yamaguti algebra, whose nontrivial brackets are given by
$$[e_1,e_2]=e_1,~\quad \Courant{e_1,e_2,e_2}=e_1,$$
where $\{e_1,e_2\}$ is a basis for $\g$.
Moreover,
$$R=
\begin{pmatrix}
0 & a\\
0 & b
\end{pmatrix}$$
is a Rota-Baxter operator on $\g$. Then by a direct computation, any element in $\wedge^2\g$ is a Nijenhuis element of $R$.
 \end{ex}

 \begin{ex}
Let $\g$ be a 4-dimensional Lie-Yamaguti algebra with a basis $\{e_1,e_2,e_3,e_4\}$ defined by
 $$[e_1,e_2]=2e_4,~\quad \Courant{e_1,e_2,e_1}=e_4.$$
 And
 $$R=
 \begin{pmatrix}
 0 & a_{12}& 0 & 0 \\
 0 & 0& 0 & 0\\
 a_{31} &a_{32} & a_{33} & a_{34}\\
 a_{41} &a_{42} & a_{43} & a_{44}
 \end{pmatrix}$$
 is a Rota-Baxter operator on $\g$. Then any element in $\wedge^2\g$ is a Nijenhuis element of $R$. In particular,
 \begin{eqnarray*}
 \frkX_1=e_1\wedge e_2,\quad \frkX_2=e_1\wedge e_3,\quad \frkX_3=e_1\wedge e_4,\\
 \frkX_4=e_2\wedge e_3,\quad \frkX_5=e_2\wedge e_4,\quad \frkX_6=e_3\wedge e_4,
 \end{eqnarray*}
 are all Nijenhuis elements of $R$.
 \end{ex}

 \subsection{Higher order deformations of relative Rota-Baxter operators}
 In this subsection, we introduce a cohomology class associated to an order $n$ deformation of a relative Rota-Baxter operator, and show that an order $n$ deformation is extendable if and only if this cohomology class is trivial. Thus we call this cohomology class the obstruction class of an order $n$ deformation being extendable.

 \begin{defi}
 Let $T$ be a relative Rota-Baxter operator on a Lie-Yamaguti algebra $(\g,[\cdot,\cdot],\Courant{\cdot,\cdot,\cdot})$ with respect to a representation $(V;\rho,\mu)$. If $T_t=\sum_{i=0}^n\frkT_it^i $ with  $\frkT_0=T$, $\frkT_i\in \Hom_{\mathbb K}(V,\g)$, $i=1,2,\cdots,n$, defines a $\mathbb K[t]/(t^{n+1})$-module from $V[t]/(t^{n+1})$ to the Lie-Yamaguti algebra $\g[t]/(t^{n+1})$ satisfying
 \begin{eqnarray*}
 [T_tu,T_tv]&=&T_t(\rho(T_t)u-\rho(T_tv)u),\\
\Courant{T_t,T_tv,T_tw}&=&T_t(D_{\rho,\mu}(T_tu,T_tv)w+\mu(T_tv,T_tw)u-\mu(T_tu,T_tw)v), \quad \forall u,v,w \in V,
 \end{eqnarray*}
 we say that $T_t$ is an {\bf order $n$ deformation} of the relative Rota-Baxter operator $T$.
 \end{defi}

 \begin{rmk}
 The left hand sides of the equations hold in the Lie-Yamaguti algebra $\g[t]/(t^{n+1})$ and the right hand sides of the equations above make sense since $T_t$ is a $\mathbb K[t]/(t^{n+1})$-module map.
 \end{rmk}

 \begin{defi}
 Let $T_t=\sum_{i=0}^n\frkT_it^i $ be an order $n$ deformation of a relative Rota-Baxter operator $T$ on a Lie-Yamaguti algebra $(\g,[\cdot,\cdot],\Courant{\cdot,\cdot,\cdot})$ with respect to a representation $(V;\rho,\mu)$. If there exists a $1$-cochain $\frkT_{n+1}\in \Hom_{\mathbb K}(V,\g)$ such that $\widetilde{T_t}=T_t+\frkT_{n+1}t^{n+1}$ is an order $n+1$ deformation of $T$, then we say that $T_t$ is {\bf extensible}.
 \end{defi}

 \begin{thm}\label{ob}
Let $T$ be a relative Rota-Baxter operator on a Lie-Yamaguti algebra $(\g,[\cdot,\cdot],\Courant{\cdot,\cdot,\cdot})$ with respect to a representation $(V;\rho,\mu)$, and $T_t=\sum_{i=0}^n\frkT_it^i$ be an order $n$ deformation of $T$. Then $T_t$ is extensible if and only if the cohomology class ~$[\Ob^T]\in \huaH_T^2(V,\g)$ is trivial, where $\Ob^T=(\Ob_I^T,\Ob_{II}^T)\in \huaC_T^2(V,\g)$ is defined to be
\begin{eqnarray*}
\Ob_I^T(v_1,v_2)&=&\sum_{i+j=n+1,\atop i,j\geqslant 1}\Big([\frkT_iv_1,\frkT_jv_2]-\frkT_i(\rho(\frkT_jv_1)v_2-\rho(\frkT_jv_2)v_1)\Big),\\
\Ob_{II}^T(v_1,v_2,v_3)&=&\sum_{i+j+k=n+1,\atop i,j,k\geqslant 1}\Big(\Courant{\frkT_iv_1,\frkT_jv_2,\frkT_kv_3}-\frkT_i(D(\frkT_jv_1,\frkT_kv_2)v_3\\
~ \nonumber&&+\mu(\frkT_jv_2,\frkT_kv_3)v_1-\mu(\frkT_jv_1,\frkT_kv_3)v_2\Big), \quad \forall v_1,v_2,v_3 \in V.
\end{eqnarray*}
\end{thm}
\begin{proof}
Let $\widetilde{T_t}=\sum_{i=0}^{n+1}\frkT_it^i$ be the extension of $T_t$, then for all $u,v,w \in V$,
\begin{eqnarray}
\label{n order1}[\widetilde{T_t}u,\widetilde{T_t}v]&=&\widetilde{T_t}\Big(\rho(\widetilde{T_t}u)v-\rho(\widetilde{T_t}v)u\Big),\\
\label{n order2}\Courant{\widetilde{T_t}u,\widetilde{T_t}v,\widetilde{T_t}w}&=&\widetilde{T_t}\Big(D(\widetilde{T_t}u,\widetilde{T_t}v)w+\mu(\widetilde{T_t}v,\widetilde{T_t}w)u
-\mu(\widetilde{T_t}u,\widetilde{T_t}w)v\Big).
\end{eqnarray}
Expanding the Eq. \eqref{n order1} and comparing the coefficients of $t^n$ yields that
\begin{eqnarray*}
\sum_{i+j=n+1,\atop i,j \geqslant 0}\Big([\frkT_iu,\frkT_jv]-\frkT_i\big(\rho(\frkT_ju)v-\rho(\frkT_jv)u\big)\Big)=0,
\end{eqnarray*}
which is equivalent to
\begin{eqnarray*}
~ &&\sum_{i+j=n+1,\atop i,j \geqslant 1}\Big([\frkT_iu,\frkT_jv]-\frkT_i\big(\rho(\frkT_ju)v-\rho(\frkT_jv)u\big)\Big)+[\frkT_{n+1}u,Tv]+[Tu,\frkT_{n+1}v]\\
~ &&\quad-\big(T(\rho(\frkT_{n+1}u)v-\rho(\frkT_{n+1}v)u)+\frkT_{n+1}(\rho(Tu)v-\rho(Tv)u)\big)=0,
\end{eqnarray*}
which is also equivalent to
\begin{eqnarray}
\Ob_I^T+\delta_I^T(\frkT_{n+1})=0.\label{ob:cocy1}
\end{eqnarray}
Similarly, expanding the Eq. \eqref{n order2} and comparing the coefficients of $t^n$ yields that
\begin{eqnarray}
\Ob_{II}^T+\delta_{II}^T(\frkT_{n+1})=0.\label{ob:cocy2}
\end{eqnarray}
From \eqref{ob:cocy1} and \eqref{ob:cocy2}, we get
$$\Ob_T=-\delta^T(\frkT_{n+1}).$$
Thus, the cohomology class $[\Ob_T]$ is trivial.

Conversely, suppose that the cohomology class $[\Ob_T]$ is trivial, then there exists $\frkT_{n+1}\in \huaC_T^1(V,\g)$, such that ~$\Ob_T=-\delta^T(\frkT_{n+1}).$ Set $\widetilde{T_t}=T_t+\frkT_{n+1}t^{n+1}$. Then for all $0 \leqslant s\leqslant n+1$, ~$\widetilde{T_t}$ satisfies
\begin{eqnarray*}
\sum_{i+j=s}\Big([\frkT_iu,\frkT_jv]-\frkT_i\big(\rho(\frkT_ju)v-\rho(\frkT_jv)u\big)\Big)=0,\\
\sum_{i+j+k=s}\Big(\Courant{\frkT_iu,\frkT_jv,\frkT_kw}-\frkT_i\big(D(\frkT_ju,\frkT_kv)w+\mu(\frkT_jv,\frkT_kw)u-\mu(\frkT_ju,\frkT_kw)v\big)\Big)=0,
\end{eqnarray*}
which implies that $\widetilde{T_t}$ is an order $n+1$ deformation of $T$. Hence it is a extension of $T_t$.
\end{proof}

\begin{defi}
Let $T$ be a relative Rota-Baxter operator on a Lie-Yamaguti algebra $(\g,[\cdot,\cdot],\Courant{\cdot,\cdot,\cdot})$ with respect to a representation $(V;\rho,\mu)$, and $T_t=\sum_{i=0}^n\frkT_it^i$ be an order $n$ deformation of $T$. Then the cohomology class $[\Ob^T]\in \huaH_T^2(V,\g)$ defined in Definition {\rm \ref{ob}} is called the {\bf obstruction class } of $T_t$ being extensible.
\end{defi}

\begin{cor}
Let $T$ be a relative Rota-Baxter operator on a Lie-Yamaguti algebra $(\g,[\cdot,\cdot],\Courant{\cdot,\cdot,\cdot})$ with respect to a representation $(V;\rho,\mu)$. If $\huaH_T^2(V,\g)=0$, then every $1$-cocycle in $\huaZ_T^1(V,\g)$ is the infinitesimal of some formal deformation of the relative Rota-Baxter operator $T$.
\end{cor}

\vspace{2mm}

 \noindent
 {\bf Acknowledgements:}
Qiao was partially supported by NSFC grant 11971282.


\begin{thebibliography}{a}















\bibitem{Bai-Bellier-Guo-Ni}
C. Bai, O. Bellier, L. Guo, and X. Ni, Splitting of operations, Manin products, and Rota-Baxter operators, \emph{Int. Math. Res. Not. IMRN}  (2013),  no. 3, 485-524.




\bibitem{bala}
D. Balavoine, Deformations of algebras over a quadratic operad. Operad: Proceedings of Renaissance Conferences (Hartford, CT/Luminy, 1995), \emph{Contemp. Math. Amer. Math. Soc.}, Providence, vol.{\bf 202}, (1997) 207-234.


\bibitem{Ba} G. Baxter, An analytic problem whose solution follows from a simple algebraic identity, \emph{Pacific J. Math.}  {\bf 10}  (1960), 731-742.

\bibitem{B.B.M}
P. Benito, M. Bremner, and S. Madariaga, Symmetric matrices, orthogonal Lie algebras and Lie-Yamaguti algebras, \emph{Linear Multilinear Algebra} {\bf 63} (2015)
1257-1287.

\bibitem{B.D.E}
P. Benito, C. Draper, and A. Elduque, Lie-Yamaguti algebras related to $\huaG_2$, \emph{J. Pure Appl. Algebra} {\bf 202} (2005) 22-54.

\bibitem{B.E.M1}
P. Benito, A. Elduque, and F. Mart$\acute{i}$n-Herce, Irreducible Lie-Yamaguti algebras, \emph{J. Pure Appl. Algebra} {\bf 213} (2009) 795-808.

\bibitem{B.E.M2}
P. Benito, A. Elduque, and F. Mart$\acute{i}$n-Herce, Irreducible Lie-Yamaguti algebras of generic type, \emph{J. Pure Appl. Algebra} {\bf 215} (2011) 108-130.



\bibitem{Burde}
D. Burde, Simple left-symmetric algebras with solvable Lie algebra, \emph{Manuscr. Math.} {\bf 95} (1988) 397-411.

\bibitem{CP}
V. Chari and A. Pressley, A Guide to Quantum Groups, Cambridge University Press, (1994) Cambridge.






\bibitem{Gerstenhaber1}
M. Gerstenhaber, The cohomology structure of an associative ring, \emph{Ann. Math.} {\bf 78}, (1963) 267-288.

\bibitem{Gerstenhaber2}
M. Gerstenhaber, On the deformations of rings and algebras, \emph{Ann. Math.}(2) {\bf 79}, (1964) 59-103.

\bibitem{Gerstenhaber3}
M. Gerstenhaber, On the deformations of rings and algebras II, \emph{Ann. Math.} {\bf 84}, (1966) 1-19.

\bibitem{Gerstenhaber4}
M. Gerstenhaber, On the deformations of rings and algebras III, \emph{Ann. Math.} {\bf 88}, (1968) 1-34.


\bibitem{Gerstenhaber5}
M. Gerstenhaber, On the deformations of rings and algebras IV, \emph{Ann. Math.} {\bf 99}, (1974) 257-276.

\bibitem{Get}
E. Getzler, Lie theory for nilpotent $L_\infty$-algebras, \emph{Ann. Math. (2)}, {\bf 170} (2009), 271-301.

\bibitem{GLST}
A. Guan, A. Lazarev, Y. Sheng, and R. Tang, Review of deformation theory I: concrete formulas for deformations of algebraic structures, \emph{Adv. Math. (China)} {\bf49} (2020) no.3, 257-277.

\bibitem{Gub}
L. Guo,  An introduction to Rota-Baxter algebra. Surveys of Modern Mathematics, 4. International Press, Somerville, MA; Higher Education Press, Beijing, 2012. xii+226 pp.

\bibitem{Hart}
R. Hartshore, Deformation Theory, \emph{Graduate Texts in Math.} 257 (2010) Springer, Berlin.



\bibitem{Weinstein}
M. K. Kinyon and A. Weinstein, Leibniz algebras, Courant algebroids and multiplications on reductive homogeneous spaces, \emph{Amer. J. Math.} {\bf 123} (2001)
525-550.

\bibitem{Kodaira}
K. Kodaira and D. Spencer, On deformations of complex analytic structures I and II, \emph{Ann. Math.} {\bf 67}, 328-466(1958).

\bibitem{Kont1}
M. Kontsevich, Operads and motives in deformation quantization, \emph{Lett. Math. Phys.} {\bf 48}, 35-72(1999).

\bibitem{Kont2}
M. Kontsevich, Deformation quantization of Poisson manifolds, \emph{Lett. Math. Phys.} {\bf 66}, 157-216(2003).

\bibitem{Kupershmidt}
B. A. Kupershmidt, What a classical $r$-matrix really is, \emph{J. Nonlinear Math. Phys.} {\bf 6} (1999) no.4, 448-488.

\bibitem{Lada2}
T. Lada and M. Markl, Strongly homotopy Lie algebras, \emph{Comm. Algebra} {\bf 23} (1995), 2147-2161.

\bibitem{Lada1}
T. Lada and J. Stasheff, Introduction to sh Lie algebras for physicists, \emph{Internat. J. Theoret. Phys.} {\bf 32} (1993), 1087-1103.

\bibitem{L.CHEN}
J. Lin, L. Chen, and Y. Ma, On the deformaions of Lie-Yamaguti algebras, \emph{Acta. Math. Sin. (Engl. Ser.)} {\bf 31} (2015) 938-946.


\bibitem{Loday}
J.-L. Loday and B. Vallette, Algebraic Operads, Springer, 2012.

\bibitem{Lister}
G. Lister, A structure theory of Lie triple systems, \emph{Trans. Amer. Math. Soc.} {\bf 72} (1952) 217-242.







\bibitem{Nij1}
A. Nijenhuis and R. Richardson, Cohomology and deformations in graded Lie algebras, \emph{Bull. Am. Math. Soc.} {\bf 72} (1966) 1-29.

\bibitem{Nij2}
A. Nijenhuis and R. Richardson, Commutative algebra cohomology and deformations of Lie and associative algebras, \emph{J. Algebra} {\bf 9} (1968) 42-105.

\bibitem{Nomizu}
K. Nomizu, Invariant affine connections on homogeneous spaces, \emph{Amer. J. Math.} {\bf76} (1954) 33-65.

\bibitem{PBG}
J. Pei, C. Bai, and L. Guo, Splitting of operads and Rota-Baxter operators on operads,
\emph{Appl. Categ. Structures}  {\bf25}  (2017),  no. 4, 505-538.

\bibitem{Rot}
M. Rotkiewicz, Cohomology ring of $n$-Lie algebras, \emph{Extracta Math.} {\bf 20} (2005), 219-232.

\bibitem{STS}
 M. A. Semonov-Tian-Shansky, What is a classical R-matrix? \emph{Funct. Anal. Appl.} {\bf 17} (1983), 259-272.




\bibitem{Sheng Zhao}
Y. Sheng, J. Zhao, and Y. Zhou, Nijenhuis operators, product structures and complex structures on Lie-Yamaguti algebras, \emph{J. Algebra Appl.} {\bf 20} (2021), no. 8, Paper No.2150146, 22 pp.

\bibitem{SZ1}
Y. Sheng and J. Zhao, Relative Rota-Baxter operators and symplectic structures on Lie-Yamaguti algebras, \emph{Comm. Algebra} {\bf50} (2022) no. 9, 4056-4073.

\bibitem{Stasheff}
J. Stasheff, Differential graded Lie algebras, quasi-Hopf algebras and higher homotopy algebras. \emph{Quantum groups (Leningrad, 1990)}, 120-137, Lecture notes in Math., 1510, \emph{Springer, Berlin}, 1992.

\bibitem{Takahashi}
N. Takahashi, Modules over quadratic spaces and representations of Lie-Yamaguti algebras, \emph{J. Lie Theory} {\bf 31} (2021) no.4, 897-932.

\bibitem{TBGS}
R. Tang, C. Bai, L. Guo, and Y. Sheng, Deformations and their controlling cohomologies of $\huaO$-operators, \emph{Comm. Math. Phys.} {\bf 368} (2019), 655-700.

\bibitem{THS}
R. Tang, S. Hou, and Y. Sheng, Lie 3-algebras and deformations of relative Rota-Baxter operators on 3-Lie algebras,
\emph{J. Algebra}  {\bf567}  (2021), 37-62.

\bibitem{T.S2}
R. Tang, Y. Sheng, and Y. Zhou, Deformations of relative Rota-Baxter operators on Leibniz algebras, \emph{Int. J. Geom. Methods Mod. Phys.} {\bf 17} (2020), no. 12, 2050174, 21 pp.

\bibitem{Voronov}
Th. Vornovo, Higher derived brackers and homotopy algebras, \emph{J. Pure Appl. Algebra} {\bf 202} (2005), 133-153.

\bibitem{Yamaguti1}
K. Yamaguti, On the Lie triple system and its generalization, \emph{J. Sci. Hiroshima Univ. Ser. A} {\bf 21} (1957/1958) 155-160.

\bibitem{Yamaguti2}
K. Yamaguti, On cohmology groups of general Lie triple systems, \emph{Kumamoto J. Sci. A} {\bf 8} (1967/1969) 135-146.

\bibitem{Zhang1}
T. Zhang and J. Li, Deformations and extension of Lie-Yamaguti algebras, \emph{Linear Multilinear Algebra.} {\bf 63} (2015) 2212-2231.

\bibitem{ZQ1}
J. Zhao and Y. Qiao, The classical Lie-Yamaguti Yang-Baxter equation and Lie-Yamaguti bialgebras (in Chinese), \emph{Sci. Sin. Math.} {\bf 49}  (2023) 1-23, doi: 10.1360/SCM-2022-0517.

\bibitem{ZQ2}
J. Zhao and Y. Qiao, Cohomology and relative Rota-Baxter-Nijenhuis structures on $\mathsf{LieYRep}$ pairs, \emph{J. Geom. Phys.} {\bf 186} (2023) No. 104749, 24pp.






\end{thebibliography}
\end{document}